 \documentstyle[12pt]{article}

\if@twoside  m
\else
\fi
\Roman{section}

\font\sc = eusm10 scaled \magstep 1
\font\gt = eufm10 scaled \magstep 1
\font\Lgt = eufm10 scaled 2488  
\font\lgt = eufm10 scaled 2074  
\font\sr = msbm10 scaled \magstep 1

\def\Uh{${\cal U}_h(${\gt su}$(N))\ $}
\def\Uha{{\cal U}_h(${\gt su}$(N))}
\def\LUh{${\cal U}_h(${\Lgt su}$(N))\ $}
\def\lUh{${\cal U}_h(${\lgt su}$(N))\ $}
\def\SU{$SU(N)\ $}
\def\Fh{${\cal F}_{hol}\ $}
\def\Fah{${\cal F}_{ahol}\ $}
\def\Faha{{\cal F}_{ahol}}
\def\Fha{{\cal F}_{hol}}
\def\Gmh{${\cal G}^{(m)}_{hol}\ $}
\def\Gmah{${\cal G}^{(m)}_{ahol}\ $}

\def\Gnh{${\cal G}^{(n)}_{hol}\ $}

\def\Fah{${\cal F}_{ahol}\ $}
\def\AqSU{${\cal A}_q(SU(N))\ $}
\def\AqAN{${\cal A}_q(AN)\ $}
\def\Ada{\tilde{\cal A}_d}
\def\Ad{$\tilde{\cal A}_d\ $}
\def\DeltaYa{\Delta Y=Y_{(1)}\otimes Y_{(2)}}
\def\DeltaY{$\Delta Y=Y_{(1)}\otimes Y_{(2)}\ $}
\def\DeltaYa{\Delta Y=Y_{(1)}\otimes Y_{(2)}}

\def\DeltaX{$\Delta X=X_{(1)}\otimes X_{(2)}\ $}

\def\Deltaf{$\Delta f=f_{(1)}\otimes f_{(2)}\ $}
\def\bF{\bar{\cal F}}
\def\sXm{\mbox{\sc X}^{(m)}}
\def\sXbm{\bar{\mbox{\sc X}}^{(m)}}
\def\sXn{\mbox{\sc X}^{(n)}}

\def\sYm{\mbox{\sc Y}^{(m)}}

\def\sYn{\mbox{\sc Y}^{(n)}}

\def\sX{\mbox{\sc X}}
\def\sXb{\bar{\mbox{\sc X}}}
\def\sY{\mbox{\sc Y}}
\def\gZm{\mbox{\gt Z}^{(m)}}
\def\bgZm{\bar{\mbox{\gt Z}}^{(m)}}
\def\gZn{\mbox{\gt Z}^{(n)}}
\def\gZ{\mbox{\gt Z}}
\def\gZmast{\mbox{\gt Z}^{(m)\raisebox{1pt}{\footnotesize$\ast$}}}
\def\gZnast{\mbox{\gt Z}^{(n)\raisebox{1pt}{\footnotesize$\ast$}}}
\def\gZast{\mbox{\gt Z}^\ast}
\def\gQm{\mbox{\gt Q}^{(m)}}
\def\gQn{\mbox{\gt Q}^{(n)}}
\def\gQmm{\mbox{\gt Q}^{(m-1)}}
\def\gQnm{\mbox{\gt Q}^{(n-1)}}
\def\gQ{\mbox{\gt Q}}
\def\gXm{\mbox{\gt X}^{(m)}}
\def\gXn{\mbox{\gt X}^{(n)}}
\def\bgXm{\bar{\mbox{\gt X}}^{(m)}}
\def\bgXn{\bar{\mbox{\gt X}}^{(n)}}
\def\Xm{X^{(m)}}
\def\Xn{X^{(n)}}
\def\bXm{\bar{X}^{(m)}}
\def\bXn{\bar{X}^{(n)}}
\def\Em{E^{(m)}}
\def\En{E^{(n)}}
\def\Fm{F^{(m)}}
\def\Fn{F^{(n)}}
\def\Fmm{F^{(m-1)}}
\def\Fnm{F^{(n-1)}}
\def\Zm{Z^{(m)}}
\def\Zn{Z^{(n)}}
\def\Am{A^{(m)}}
\def\Bm{B^{(m)}}
\def\Cm{C^{(m)}}
\def\Dm{D^{(m)}}
\def\Zmast{Z^{(m)\raisebox{1pt}{\footnotesize$\ast$}}}

\def\zm{z^{(m)}}
\def\zmast{z^{(m)\raisebox{1pt}{\footnotesize$\ast$}}}
\def\zsast{z^{(s)\raisebox{1pt}{\footnotesize$\ast$}}}
\def\bzm{\bar{z}^{(m)}}
\def\Ejj{E_{jj}}
\def\Ejjp{E_{j,j+1}}
\def\Ejpj{E_{j+1,j}}
\def\Ejpjp{E_{j+1,j+1}}
\def\qHj{q^{H_{j}}}
\def\qmHj{q^{-H_{j}}}
\def\qpmHj{q^{\pm H_{j}}}
\def\Xpj{X^+_j}
\def\Xmj{X^-_j}
\def\Xpmj{X^\pm_j}
\def\qHk{q^{H_{k}}}

\def\Xmk{X^-_k}
\def\Xpmk{X^\pm_k}
\def\Aqa{{\cal A}_q}
\def\cR{{\cal R}}

\def\Ca{\mbox{\sr C}}

\def\Ra{\mbox{\sr R}}
\def\gZa{\mbox{\gt Z}}

\def\sN{{\sc N}\ }
\def\sNa{\mbox{\sc N}}
\def\sZ{{\sc Z}\ }
\def\sZa{\mbox{\sc Z}}
\def\cF{$\cal F\ $}
\def\cFa{{\cal F}}
\def\e{\mbox{e}}
\def\vare{\varepsilon}
\def\bI{{\bf I}}
\def\bO{{\bf 0}}
\def\span{\mbox{span}}
\def\diag{\mbox{diag}}

\def\fsgn{\footnotesize\mbox{sgn}}
\def\Mat{\mbox{Mat}}
\def\LHS{\mbox{LHS}}
\def\RHS{\mbox{RHS}}
\def\id{\mbox{id}}
\def\tr{\mbox{tr}}
\def\End{\mbox{End}}
\newcommand{\QED}{\mbox{\rule[-1.5pt]{6pt}{10pt}}}

\def\Tla{{\cal T}^\lambda}
\def\Ml{${\cal M}_\lambda\ $}
\def\Mla{{\cal M}_\lambda}
\def\Hl{${\cal H}_\lambda\ $}
\def\Hla{{\cal H}_\lambda}
\def\wl{$w_\lambda\ $}
\def\wla{w_\lambda}
\def\el{$e_\lambda\ $}
\def\ela{e_\lambda}
\def\tl{$\tau_\lambda\ $}

\begin{document}

\title{\bf Representations of \LUh derived from quantum flag manifolds}
\author{Pavel \v S\v tov\'\i\v cek\thanks{Humboldt Fellow}
\thanks{On leave of absence from: Department of Mathematics,
Faculty of Nuclear Science, CTU, Prague, Czech Republic}\
\ and Reidun Twarock}
\date{}
\maketitle

\begin{quote}
\begin{center}
{\small Arnold Sommerfeld Institute\\
TU Clausthal, Leibnizstr. 10\\
D-38678 Clausthal-Zellerfeld, Germany}
\end{center}
\end{quote}
\begin{abstract}
\noindent
{\small A relationship between quantum flag and Grassmann manifolds
is revealed. This enables a formal diagonalization of quantum positive
matrices. The requirement that this diagonalization defines a homomorphism
leads to a left \Uh -- module structure on the algebra generated by
quantum antiholomorphic coordinate functions living on the flag
manifold. The module is defined by prescribing the action on the unit
and then extending it to all polynomials using a quantum version of
Leibniz rule. Leibniz rule is shown to be induced by the dressing
transformation. For discrete values of parameters occuring in the
diagonalization one can extract finite-dimensional irreducible
representations of \Uh as cyclic submodules.}
\end{abstract}

\section{INTRODUCTION}

Flag manifolds were quantized already some time
ago \cite{T-T, Soibelman}. Also some other types of homogeneous spaces
were treated including Grassmann manifolds \cite{Podlesz, Sheu, JMP, QG}.
Moreover, the quantization can be described in a unified way for all types
of coadjoint orbits regarded as complex manifolds and for all simple compact
groups from the four principal series \cite{J-S-II}. Quantum
homogeneous spaces were related to representations and co-representations
when taking various points of view like that of induced representations,
utilizing $q$-difference operators
etc. \cite{Parshall-W, Japanes, Biedenharn-L, Fiore}.
Other aspects are of interest, too, like applications to physical models and
differential geometry \cite{Brzezinski-M, Schuppetal, Woronowicz, Jurco}.
In fact, we just succeeded to quote only a small part of contributions
related to this subject (cf. Ref. \cite{Chari-P}).

In a recent paper \cite{Jurco-Sch}, carrier spaces of representations of
\Uh were realized as subspaces in the algebra generated by quantum
coordinate functions on the flag manifold. The particular case of
$SU(3)$ was treated in an analogous way in Ref. \cite{Bulgariens}.
Our goal is to derive a similar description but taking a different
approach and presenting some additional results, too. Also the obtained
formulas are optically somewhat different though necessarily convertible
one into each other. While the method of Ref. \cite{Jurco-Sch} is based
 on a decomposition of the universal R-matrix
\cite{Kirillov-R, Levendorskii-S, Khoros-Tolstoy} we start directly
from the quantized flag manifold.

Before explaining main features of
our approach let us devote a couple of words to the classical case.
Quantum flag manifolds can be viewed as quantized orbits of the classical
dressing transformation of \SU acting on its solvable dual $AN$
\cite{SemenovTS, Lu-W}. The solvable group $AN$ comes from the Iwasawa
decomposition $SL(N,\Ca)=SU(N)\cdot AN$ and is formed by unimodular
upper-triangular matrices having positive diagonals. To make the structure
of dressing orbits more transparent one can use the bijection sending
$\Lambda\in AN$ to a unimodular positive matrix
$M=\Lambda^\ast\Lambda$. The dressing
action on positive matrices becomes just the unitary transformation
$M\mapsto U^\ast MU$ and thus the orbits are determined by sets of
eigen-values. We attempted to find a quantum analog to this procedure.

To diagonalize quantum positive matrices we needed, first of all, to reveal
a relationship between quantum flag and Grassmann manifolds. It is quite
straightforward to see that the quantum Grassmannians are embedded into
and jointly generate the quantum algebra related to the flag manifold.
Also the constraints reducing the number of generators obtained this
way are easy to find. What is less obvious are cross commutation relations
between different Grassmannians. The decomposition of the quantum flag
manifold into Grassmann manifolds induces a diagonalization of quantum
positive matrices provided a set of parameters (eigen-values) has been
chosen. The requirement that this diagonalization defines a homomorphism
from $\Aqa(AN)$ onto the quantum flag manifold leads to a left \Uh --
module structure on the subalgebra of
the latter algebra generated by quantum
``antiholomorphic'' functions. The module is defined by prescribing the
action on the unit and then extending it to all polynomials in
non-commutative variables (quantum antiholomorphic coordinate functions)
using a recursive rule, an idea utilized already in
Ref. \cite{JALG}. Moreover, we prove that this recursive rule follows
from the quantum dressing transformation, making the role of the dressing
transformation quite explicit. Up to this point, we employed the
Faddeev--Reshetikhin--Takhtajan (FRT) description of deformed enveloping
algebras \cite{F-R-T}. However this result enables us to transcribe,
quite straightforwardly, all expressions in terms of Chevalley generators,
too. Let us note that in Ref. \cite{Jurco-Sch} only the FRT picture has
been treated. Finite-dimensional irreducible representations of \Uh are
then easily obtained as cyclic submodules with unit as the cyclic vector
and, at the same time, the lowest weight vector.

We have just explained the basic ideas following the structure of the
paper. Let us summarize that the notation is introduced and some
preliminary facts are reviewed in Section 2, the relationship between
quantum flag and Grassmann manifolds is described in Section 3, the left
\Uh -- module structure on the quantum flag manifold is derived in Section 4,
the role of the quantum dressing transformation is revealed in Section 5
and Section 6 is devoted to the transcription of all formulas in terms
of Chevalley generators as well as to finite-dimensional irreducible
representations of \Uh.

\section{PRELIMINARIES, NOTATION}

Let us recall some basic and well known facts related to the quantum
group $SU_q(N)$ and the deformed enveloping algebra \Uh
\cite{Drinfeld, Jimbo1, Jimbo2} introducing this way the notation. The
deformation parameter is $q=\e^{-h}$, with $h\in\Ra$ (or one can
consider, too, $h$ as a formal variable but real like, i.e.,
$h^\ast=h$ and to work with the ring $\Ca[[h]]$).

The $\ast$-Hopf algebra
of quantum functions living on \SU is denoted by \AqSU and $U$ stands
for the defining $N\times N$ vector corepresentation. The symbols
$\vare$ and $\Delta$ designate everywhere the counit and the
comultiplication, respectively, for any Hopf algebra under consideration
and we use Sweedler notation (a summation is understood)
$$
\DeltaYa\,.
$$
Thus
\begin{eqnarray*}
\makebox[7em]{}
& R_{12}U_1U_2=U_2U_1R_{12},\quad {\det}_qU=1,\quad U^\ast=U^{-1}, &
\makebox[6em][r]{}\cr
\makebox[7em]{}
& \vare(U)=\bI,\quad \Delta U=U\dot{\otimes}U\,. &
\,\makebox[6em][r]{(2.1)}\cr
\end{eqnarray*}
Here we define, as usual,
$(A\dot{\otimes}B)_{jk}:=\sum_s A_{js}\otimes B_{sk}$ and
$(A^\ast)_{jk}:=(A_{kj})^\ast$.

The R-matrix acting in $\Ca^N\otimes\Ca^N$
and obeying Yang--Baxter (YB) equation
$R_{12}R_{13}R_{23}=R_{23}R_{13}R_{12}$ is given by
$$
R_{jk,st}:=\delta_{js}\delta_{kt}+(q-q^{\fsgn(k-j)})
\delta_{jt}\delta_{ks}\,.
\eqno(2.2)$$
In what follows we identify linear operators in $\Ca^N$ (and its tensor
products) with their matrices in the standard basis $\{e_1,e_2,\dots,e_N\}$.
In tensor products the lexicographical ordering is assumed. The R-matrix
verifies the relations
$$
R_{q^{-1}}=R_q^{-1},\quad R_{12}^t= R_{21}\equiv PR_{12}P
\eqno(2.3)$$
and Hecke condition
$$
(q-q^{-1})P=R_{12}-R_{21}^{-1}=R_{21}-R_{12}^{-1}\,.
\eqno(2.4)$$
Here $P$, $P_{jk,st}:=\delta_{jt}\delta_{ks}$, is the flip operator.

The $\ast$-Hopf algebra \Uh is generated by Chevalley generators
$\qpmHj,\ \Xpj,\ \Xmj$, $1\le j\le N-1$, and is determined by the
relations
\begin{eqnarray*}
&& \makebox[3em][r]{}
[\qHj,\qHk]=0\,,
\makebox[9em][r]{}\cr
&& \makebox[3em][r]{}
 \qHj\Xpmk=
 \left\{\begin{array}{c} q^{\pm2} \\ q^{\mp1} \\ 1 \end{array}\right\}
\times\Xpmk\qHj , \quad
for\ \left\{\begin{array}{l} j=k\ ,\\ |j-k|=1, \\ |j-k|\ge2,\end{array}\right.
\makebox[9em][r]{(2.5)}\cr
&& \makebox[3em][r]{}
 [\Xpj,\Xmk]=\delta_{jk}(\qHj-\qmHj)/(q-q^{-1})\,,
\makebox[9em][r]{}\cr
&& \makebox[3em][r]{}
(\Xpmj)^2\Xpmk-(q+q^{-1})\Xpmj\Xpmk\Xpmj+\Xpmk(\Xpmj)^2=0,\quad
for\ |j-k|=1,
\makebox[9em][r]{}\cr
&& \makebox[3em][r]{}
[\Xpmj,\Xpmk]=0,\quad for\ |j-k|\ge2,
\makebox[9em][r]{}\cr
\end{eqnarray*}
and
$$
(\qpmHj)^\ast=\qpmHj,\quad (\Xpmj)^\ast=X^\mp_j\,.
\eqno(2.6)$$
Furthermore,
\begin{eqnarray*}
\makebox[3em]{}
& \vare(\qpmHj)=1,\quad \vare(\Xpmj)=0\,,  &
\makebox[3em][r]{(2.7)}\cr
\makebox[3em]{}
& \Delta(\qpmHj)=\qpmHj\otimes\qpmHj,\quad
\Delta(\Xpmj)=\Xpmj\otimes q^{-H_j/2}+q^{H_j/2}\otimes\Xpmj\,.  &
\makebox[3em][r]{}\cr
\end{eqnarray*}

There exists another description of \Uh due to
Faddeev--Reshetikhin--Takhtajan \cite{F-R-T} which can be reinterpreted
as the quantization of the generalized dual of \SU, namely the solvable
group $AN$ coming from the Iwasawa decomposition
$SL(N,\Ca)=SU(N)\cdot AN$. The $\ast$-Hopf algebra \AqAN is generated
by entries of the upper triangular matrix $\Lambda=(\alpha_{jk})$ and its
adjoint $\Lambda^\ast$ and is determined by the relations
\begin{eqnarray*}
\makebox[7em]{}
& R_{12}\Lambda_1\Lambda_2=\Lambda_2\Lambda_1R_{12},\quad
\Lambda_1^\ast R_{12}^{-1}\Lambda_2=\Lambda_2R_{12}^{-1}\Lambda_1^\ast\,, &
\makebox[7em][r]{}\cr
\makebox[7em]{}
& \alpha_{jj}^\ast=\alpha_{jj},\quad \prod\limits_{j=1}^N\alpha_{jj}=1\,, &
\makebox[7em][r]{(2.8)}\cr
\makebox[7em]{}
& \vare(\Lambda)=\bI,\quad\Delta_{AN}(\Lambda)=\Lambda\dot{\otimes}\Lambda\,,
&
\makebox[7em][r]{}\cr
\end{eqnarray*}
The $\ast$-algebras \Uh and \AqAN can be identified by an isomorphism which
is given explicitly in terms of generators ($1\le j\le N-1$),
\begin{eqnarray*}
\makebox[7em]{}
\qHj & = & \alpha_{jj}^{-1}\alpha_{j+1,j+1}\,,
\makebox[7em][r]{}\cr
\makebox[7em]{}
(q-q^{-1})\Xpj & = & -q^{-1/2}(\alpha_{jj}\alpha_{j+1,j+1})^{-1/2}
\alpha_{j,j+1}^\ast\,,
\makebox[7em][r]{(2.9)}\cr
\makebox[7em]{}
(q-q^{-1})\Xmj & = & -q^{-1/2}(\alpha_{jj}\alpha_{j+1,j+1})^{-1/2}
\alpha_{j,j+1}\,,
\makebox[7em][r]{}\cr
\end{eqnarray*}
We note that the diagonal elements $\alpha_{jj}$ mutually commute and
$\alpha_{jj}\alpha_{j+1,j+1}$ commutes with $\alpha_{j,j+1}$. However
\Uh and \AqAN are opposite as coalgebras. Thus one has to be careful about
the comultiplication and this is why we emphasized its origin in (2.8).
To avoid any ambiguity, throughout the text $\Delta$ is always assumed
to come from \Uh rather than from \AqAN.

It turns out as more convenient to work with the $\ast$-algebra
$\Ada\subset\Aqa(AN)$ generated by entries of the matrix
$M:=\Lambda^\ast\Lambda$ rather than directly with \AqAN. It is
straightforward to see that \Ad is determined by the relations
$$
M_2R_{12}^{-1}M_1R_{21}^{-1}=R_{12}^{-1}M_1R_{21}^{-1}M_2,\quad
M^\ast=M\,.
\eqno(2.10)$$
In fact, one is not loosing that much as the algebra \AqAN can be, in
principle, again recovered from \Ad when decomposing $M$ into a product of
lower and upper triangular matrices (cf. Proposition 3.2 in Ref. \cite{QG}).

An important feature is the duality between \Uh and \AqSU expressed in
terms of a non-degenerate pairing. The both structures are combined
according to the rules
$$
\langle X,fg\rangle=\langle X_{(1)},f\rangle\langle X_{(2)},g\rangle,\quad
\langle XY,f\rangle=\langle X,f_{(1)}\rangle\langle Y,f_{(2)}\rangle\,,
\eqno(2.11)$$
with \DeltaX, \Deltaf, and
$$
\langle X,1\rangle=\vare(X),\quad\langle 1,f\rangle=\vare(f)\,.
\eqno(2.12)$$
The pairing can be described explicitly in terms of generators. Let
$E_{jk}$ be the matrix units acting as rank-one operators:
$E_{jk}v:=(e_k,v)\,e_j$ (the indices shouldn't be confused with the leg
notation referring to tensor products). Then we have
$$
\langle\qpmHj,U\rangle=q^{\pm(\Ejj-\Ejpjp)},\
\langle\Xpj,U\rangle=\Ejjp,\ \langle\Xmj,U\rangle=\Ejpj\,.
\eqno(2.13)$$
In the FRT picture we have
$$
\langle\Lambda_1,U_2\rangle=R_{21}^{-1},\quad
\langle\Lambda^\ast_1,U_2\rangle=R_{12}^{-1}\,.
\eqno(2.14)$$

To describe a relationship between quantum flag and Grassmann manifolds
we shall need the following family of orthogonal projectors. Let $\Em$ be
the matrix corresponding to the orthogonal projector in $\Ca^N$ onto
$\span\{e_1,\dots,e_m\}$ and set $\Fm:=\bI-\Em$. Thus
$$
\Em=\sum_{j\le m}\Ejj,\quad\Fm=\sum_{j>m}\Ejj\,;
$$
particularly, $E^{(0)}=\bO,\ E^{(N)}=\bI$.
Quite important is the following relation between $\Em$ and the R-matrix,
$$
\Em_1R_{12}=\Em_1R_{12}\Em_1\,,
\eqno(2.15a)$$
and consequently,
$$
R_{12}\Em_2=\Em_2R_{12}\Em_2,\ R_{12}\Fm_1=\Fm_1R_{12}\Fm_1,\
\Fm_2R_{12}=\Fm_2R_{12}\Fm_2\,.
\eqno(2.15b)$$
Observe also that
$$
\Em_1\Fm_2R_{12}=\Em_1\Fm_2\,,
\eqno(2.16)$$
and
$$
\Em_1\Em_2R_{12}=R_{12}\Em_1\Em_2,\quad
\Fm_1\Fm_2R_{12}=R_{12}\Fm_1\Fm_2\,.
\eqno(2.17)$$

\section{QUANTUM FLAG AND GRASSMANN\newline
MANIFOLDS}

There is a standard way of introducing local holomorphic coordinate functions
on homogeneous spaces $SU(N)/S(U(m_1)\times\dots\times U(m_k))$,
$\sum m_j=N$, which is given by Gauss decomposition.
This coordinate system is well defined on the so called big cell (the unique
cell of the top dimension in the cell decomposition) which covers the whole
manifold up to an algebraic subset. The coordinate functions appear as
entries of a block upper-triangular matrix \sN with unit blocks on the
diagonal. The structure of the blocks depends on the type of the homogeneous
space. The quantization procedure for the algebra of (anti)holomorphic
functions living on the big cell has been performed successfully in many
particular cases \cite{T-T, Soibelman, Japanes, QG}. But there is a unified
and compact way of writing down the commutation relations which is valid
for any homogeneous space of the above type, namely \cite{J-S-II}
$$
R_{12}Q_{12}^{-1}\sNa_1Q_{12}\sNa_2=Q_{21}^{-1}\sNa_2Q_{21}\sNa_1R_{12}\,,
\eqno(3.1)$$
where the matrix $Q$ is obtained from $R$ by annulating some entries in
dependence on the concrete homogeneous space in question.
Let us specialize (3.1) to flag and Grassmann manifolds.

In the case of the flag manifold, \sN is simply an upper triangular matrix
with units on the diagonal and we redenote it as \sZ,
$$
\sZa=(\zeta_{jk}),\ 1\le j,k\le N,\ with\ \zeta_{jj}=1,\ \zeta_{jk}=0\
for\ j>k\,,
\eqno(3.2)$$
and $Q=\diag(R),\ \diag(R)_{jk,st}=q^{\delta_{jk}}\delta_{js}\delta_{kt}$.
As $\diag(R)$ commutes with $R$ it is possible to simplify (3.1),
$$
R_{12}\sZa_1\diag(R)\sZa_2=\sZa_2\diag(R)\sZa_1R_{12}\,.
\eqno(3.3)$$
One can rewrite (3.3) in terms of matrix entries,
$$
q^{\delta_{ks}}\zeta_{js}\zeta_{kt}-q^{\delta_{jt}}\zeta_{kt}\zeta_{js}=
(q^{\fsgn(k-j)}-q^{\fsgn(s-t)})q^{\delta_{js}}\zeta_{ks}\zeta_{jt}\,.
\eqno(3.4)$$
The relations (3.3) (or (3.4)) define an algebra of quantum holomorphic
functions generated by $\zeta_{jk},\ 1\le j<k\le N$, and denoted here by \Fh
while the adjoint relations define an algebra of quantum antiholomorphic
functions generated by $\zeta^\ast_{jk},\ 1\le j<k\le N$, and denoted
by \Fah.

In the case of the Grassmann manifold formed by $m$-dimensional subspaces in
$\Ca^N,\ 1\le m\le N-1$, \sN has the block structure given by
$$
\sNa=\left(\begin{array}{cc} \bI & \Zm \\
                             \bO & \bI \end{array}\right)\,,\
with\ \Zm=(\zm_{jk}),\ 1\le j\le m,\ m+1\le k\le N\,.
\eqno(3.5)$$
Thus the dimension of the block $\Zm$ is $m\times(N-m)$. Here and everywhere
in what follows we don't specify the dimensions of zero and unit blocks
for all cases are determined implicitly in an unambiguous way. We set also
$$
\gZm:=\left(\begin{array}{cc} \bO & \Zm \\
                             \bO & \bO \end{array}\right) \,,\
\eqno(3.6)$$
so that $\sNa=\bI+\gZm$ and $\gZm=\Em\gZm\Fm$. Whenever convenient we shall
define both $\gZa^{(0)}$ and $\gZa^{(N)}$ as zero matrices. The matrix $Q$
is now related to $R$ according to (cf. (2.15))
$$
Q_{12}:=\Em_2R_{12}\Em_2+\Fm_2R_{12}\Fm_2=R_{12}\Em_2+\Fm_2R_{12}\,.
\eqno(3.7a)$$
It holds also true that
$$
Q_{12}=\Em_1R_{12}\Em_1+\Fm_1R_{12}\Fm_1=\Em_1R_{12}+R_{12}\Fm_1
\eqno(3.7b)$$
and $Q_q^{-1}=Q_{q^{-1}}$. To simplify (3.1) it suffices to multiply this
relation by $\Em_1\Em_2$ from the left, to use (2.15), (2.17), (3.7) and to
observe that $\Em\sNa=\Em+\gZm$. The result is
$$
R_{21}(\Em_1+\gZm_1)R_{12}(\Em_2+\gZm_2)=(\Em_2+\gZm_2)R_{21}(\Em_1+\gZm_1)
R_{12}\,.
\eqno(3.8)$$
In fact, this relation is equivalent to \cite{QG}
$$
R_{21}^{[m]}\Zm_1\Zm_2=\Zm_2\Zm_1R_{12}^{[N-m]}\,,
\eqno(3.9)$$
where $R^{[m]}$ stands for the $m^2\times m^2$ R--matrix related to the
quantum group $SU_q(m)$ and $R^{[1]}:=q$ . One can rewrite (3.8) in terms
of entries $\zm_{jk}$,
$$
\zm_{jk}\zm_{st}-\zm_{st}\zm_{jk}=
(q^{\fsgn(j-s)}-q^{\fsgn(k-t)})\zm_{sk}\zm_{jt}\,.
\eqno(3.10)$$
For a given $m,\ 1\le m\le N-1$, the relations (3.8) (or (3.9) or (3.10))
define an algebra of quantum holomorphic functions generated by
$\zm_{jk},\ 1\le j\le m<k\le N$, and denoted by \Gmh while the adjoint
relations define an algebra of quantum antiholomorphic functions generated
by ${\zm_{jk}}^\ast,\ 1\le j\le m<k\le N$, and denoted by \Gmah.

The main goal of this section is to express \Fh in terms of \Gmh,
$1\le m\le N-1$ (and analogously for the antiholomorphic versions).
It is quite straightforward to embed the algebras \Gmh into \Fh as well as
to find the constraining relations (the algebras \Gmh are not mutually
independent as subalgebras in \Fh). A more difficult problem is to
determine the cross commutation relations between \Gmh and \Gnh for
$m\not=n$. We introduce some additional notation used only locally in this
section and only for the sake of derivation of these relations. Set
$$
\sXm:=\sZa\Em=\Em\sZa\Em,\ \sXbm:=\sZa^{-1}\Em=\Em\sZa^{-1}\Em,\
\sYm:=\Em\sZa\Fm\,.
\eqno(3.11)$$
Thus we have $\Em\sZa=\sXm+\sYm$ and
$$
\sZa^{-1}\Em\sZa=\sXbm(\sXm+\sYm),\quad\Em=\sXbm\sXm=\sXm\sXbm\,.
\eqno(3.12)$$
Next we will rewrite the commutation relation (3.3)
in terms of $\sXm$ and $\sYm$.
\proclaim Lemma 3.1.
Assume that $0\le m\le n\le N$. It holds true that
\begin{eqnarray*}
\makebox[1em]{}
R_{21}^{-1}\sXm_1\diag(R)\sXn_2 & = & \sXn_2\diag(R)\sXm_1R_{21}^{-1}\Em_1\,,
\makebox[11em][r]{(3.13)}\cr
\makebox[1em]{}
R_{21}^{-1}\sXm_1\diag(R)\sYn_2 & = & \sYn_2\sXm_1\,,
\,\,\makebox[18em][r]{(3.14)}\cr
\makebox[1em]{}
\Em_1R_{12}\sYm_1\diag(R)\sXn_2 & = & \sXn_2\diag(R)\sYm_1R_{12}\Fm_1\En_2\,,
\,\,\,\makebox[9em][r]{(3.15)}\cr
\makebox[1em]{}
\Em_1R_{12}\sYm_1\diag(R)\sYn_2 & = & \sYn_2\sYm_1R_{12}+
(q-q^{-1})\sXn_2\Fm_2\sYm_1\Fn_1P\,.
\!\!\makebox[4em][r]{(3.16)}\cr
\end{eqnarray*}

\noindent {\em Proof.}
To derive (3.13) -- (3.16) it suffices to multiply the equation (3.3) by
$\Em_1\En_2$ from the left in all cases and by $\Em_1\En_2$,
$\Em_1\Fn_2$, $\Fm_1\En_2$, $\Fm_1\Fn_2$, respectively, from the right.
One has to employ (2.15), (2.16), (2.17) and, where convenient, (2.4) to
commute the projectors with the R-matrix getting consequently, for example,
$$
\Em_1\En_2R_{12}=\Em_1\En_1\En_2R_{12}=\Em_1R_{12}\Em_1\En_2
$$
(in all cases) and similar relations like
$$
R_{12}\Em_1\Fn_2=\Em_1\Fn_2+(q-q^{-1})P\Em_1\Fn_2
$$
(in the case of (3.14)). Note also that $\sXm=\sXn\Em$ etc. We omit
further details.
\quad\QED

We shall need also the following lemma whose verification is quite easy.
\proclaim Lemma 3.2.
It holds true that
$$
\En_1R_{12}^{\pm1}F^{(n-1)}_1=\diag(R)^{\pm1}\En_1F^{(n-1)}_1,\
F^{(n-1)}_2R_{12}^{\pm1}\En_2=\diag(R)^{\pm1}\En_2F^{(n-1)}_2\,,
\eqno(3.17)$$
and consequently,
$$
m\ge n\Longrightarrow R_{12}^{\pm1}\Em_1F^{(m-1)}_1\En_2F^{(n-1)}_2=
\diag(R)^{\pm1}\Em_1F^{(m-1)}_1\En_2F^{(n-1)}_2\,.
\eqno(3.18)$$
Evidently,
$$
m\not=n\Longrightarrow \diag(R)^{\pm1}\Em_1F^{(m-1)}_1\En_2F^{(n-1)}_2=
\Em_1F^{(m-1)}_1\En_2F^{(n-1)}_2\,.
\eqno(3.19)$$

\proclaim Proposition 3.3.
The algebras \Gmh, $1\le m\le N-1$, are embedded into \Fh via the equalities
(abusing notation, the LHS should be understood as the result of the
embedding)
$$
\Em+\gZm=\sZa^{-1}\Em\sZa\,.
\eqno(3.20)$$
Moreover, being embedded into \Fh, the subalgebras \Gmh generate jointly
\Fh for one can express
$$
\sZa=\sum_{m=1}^N F^{(m-1)}(\Em+\gZm)=
\bI+\sum_{m=1}^{N-1}F^{(m-1)}\gZm\,.
\eqno(3.21)$$
Thus, on the other hand, \Fh is isomorphic to the algebra generated by the
entries of the blocks $\Zm,\ 1\le m\le N-1$, and determined by the relations
\begin{eqnarray*}
\makebox[0em]{}
&& 1\le m<n\le N-1  \Longrightarrow
(\Em+\gZm)(\En+\gZn)= \Em+\gZm\,,
\!\!\makebox[5em][r]{(3.22)}\cr
\makebox[0em]{}
&& 1\le m\le n\le N-1 \Longrightarrow
\makebox[4em][r]{}\cr
\makebox[0em][r]{}
&& \quad R_{21}(\Em_1+\gZm_1)R_{12}(\En_2+\gZn_2)=(\En_2+\gZn_2)R_{21}
(\Em_1+\gZm_1)R_{12}\,.
\makebox[3em][r]{(3.23)}\cr
\end{eqnarray*}

{\em Remarks}. (1) The equality (3.20) makes sense also for $m=0$ and
$m=N$ when it reduces to the trivial identities $\bO=\bO$ and $\bI=\bI$,
respectively. A similar remark applies also for other equalities.

(2) Note that
$$
\sZa^{-1}\Em\sZa=\Em\sZa^{-1}\Em\sZa\quad and\quad
\sZa^{-1}\Em\sZa\Em=\sZa^{-1}\sZa\Em=\Em
$$
and thus the RHS of (3.20) has the same structure of blocks as the LHS.
Furthermore, the inverted relation (3.21) follows immediately from the fact
that $\diag(\sZa)=\bI$ and so
$$
\Fmm\sZa^{-1}\Em=\Fmm\Em=\Em-E^{(m-1)}\,,
$$
which implies
$$
\Fmm(\Em+\gZm)=(\Em-E^{(m-1)})\sZa\,.
$$

(3) (3.22) can be rewritten in a way coinciding formally with the classical
constraint,
$$
m<n\Longrightarrow \left(\begin{array}{cc} \bI & \Zm \end{array}\right)
\left(\begin{array}{c} -\Zn \\ \bI \end{array}\right) =0\,,
\eqno(3.24)$$
and, in fact, it reduces the number of generators. Observe also that it
holds trivially true, just by the structure of blocks (cf. (3.6)), that
$$
m\le n \Longrightarrow (\En+\gZn)(\Em+\gZm)= \Em+\gZm \,.
\eqno(3.25)$$

(4) (3.23) for $m=n$ reproduces the defining relation (3.8) of \Gmh and for
$m<n$ gives the desired cross commutation relations between different
quantum Grassmannians.

{\em Proof}. We assume that $m\le n$ and write, for simplicity of notation,
$\sX$ instead of $\sXm$, $\sX'$ instead of $\sXn$ and similarly for
$\sY,\ \sY',\ E,\ E'$ etc. Thus we have, for example, $EE'=E'E=E$.

Let us first show that the matrices $\gZm$ defined by the RHS of (3.20)
verify (3.22). Since $\sX+\sY=E(\sX'+\sY')$ and
$(\sXb(\sX+\sY))^2=\sXb(\sX+\sY)$ we have
\begin{eqnarray*}
(E+\gZ)(E'+\gZ') & = & \sXb(\sX+\sY)\sXb'(\sX'+\sY')=
\sXb E\sX'\sXb'(\sX'+\sY')\sXb'(\sX'+\sY') \cr
& = & \sXb(\sX+\sY)=E+\gZ\,.
\end{eqnarray*}

Next we will show that the same matrices $\gZm$ verify also (3.23). One can
derive successively
\begin{eqnarray*}
\diag(R)\sY_1R_{12}F_1\sXb'_2\sY'_2 & =  &
\sXb'_2E_1R_{12}\sY_1\diag(R)\sY'_2
\!\makebox[12em][r]{, by (3.15),} \nopagebreak \cr
& = & \sXb'_2\sY'_2\sY_1R_{12}+(q-q^{-1})E'_2F_2\sY_1F'_1P
\makebox[6em][r]{, by (3.16),}\cr
\end{eqnarray*}
whence ($\diag(R)^{-1}F_2E_1=F_2E_1$)
$$
\sXb_1\sY_1R_{12}F_1\sXb'_2\sY'_2=\sXb_1\diag(R)^{-1}\sXb'_2\sY'_2\sY_1R_{12}
+(q-q^{-1})\sXb_1\sY_1F'_1E'_2F_2P\,.
\eqno(3.26)$$
Furthermore,
\begin{eqnarray*}
\sXb_1\diag(R)^{-1}\sXb'_2\sY'_2\sX_1 & = &
\sXb_1\diag(R)^{-1}\sXb'_2R_{21}^{-1}\sX_1\diag(R)\sY'_2
\makebox[8em][r]{, by (3.14),}\cr
& = & \sXb_1\diag(R)^{-1}E'_2\diag(R)\sX_1R_{21}^{-1}E_1\sXb'_2\sY'_2
\!\makebox[6em][r]{, by (3.13),}\cr
& = & E_1E'_2R_{21}^{-1}E_1\sXb'_2\sY'_2
\makebox[6em][r]{}\cr
& = & R_{21}^{-1}E_1\sXb'_2\sY'_2\,,
\makebox[6em][r]{}\cr
\end{eqnarray*}
whence
$$\sXb_1\diag(R)^{-1}\sXb'_2\sY'_2=R_{21}^{-1}\sXb'_2\sY'_2\sXb_1\,.
\eqno(3.27)$$
Combination of (3.26) and (3.27) yields ($\gZ=\sXb\sY$)
$$
\gZ_1R_{12}F_1\gZ'_2=R_{21}^{-1}\gZ'_2\gZ_1R_{12}+
(q-q^{-1})\gZ_1F'_1E'_2F_2P\,.
$$
This relation can be rewritten as
$$
R_{21}\gZ_1R_{12}\gZ'_2=\gZ'_2R_{21}\gZ_1R_{12}+
(q-q^{-1})R_{21}\gZ_1E_2\gZ'_1P+
(q-q^{-1})R_{21}\gZ_1F'_1E'_2F_2P\,.
$$
To see this it suffices to notice (for $\gZ=E\gZ F$) that
$$
F'_2R_{21}E_1=E_1F'_2,\quad F_1R_{12}F_1=F_1R_{12}-(q-q^{-1})F_1E_2P\,.
$$
Next observe that (3.22) (already proven) means
$$
\gZ\gZ'=\gZ F'-E\gZ'\,.
$$
Thus we arrive at ($E=EE'$)
$$
R_{21}\gZ_1R_{12}\gZ'_2-\gZ'_2R_{21}\gZ_1R_{12}=
(q-q^{-1})(R_{21}\gZ_1F'_1E'_2P-R_{21}E_1\gZ'_1E_2P)\,.
\eqno(3.28)$$
To pass from (3.28) to (3.23) one can use that
$$
R_{21}E_1R_{12}E'_2=E'_2R_{21}E_1R_{12}E'_2=E'_2R_{21}E_1R_{12}
$$
as well as
\begin{eqnarray*}
\gZ'_2R_{21}E_1R_{12}-R_{21}E_1R_{12}\gZ'_2 & = &
-(q-q^{-1})R_{21}E_1\gZ'_1E_2P, \cr
E'_2R_{21}\gZ_1R_{12}-R_{21}\gZ_1R_{12}E'_2 & = &
(q-q^{-1})R_{21}\gZ_1F'_1E'_2P .\cr
\end{eqnarray*}
The last two relations follow respectively from
$$
R_{21}E_1=E_1R_{12}^{-1}+(q-q^{-1})E_1E_2P,\quad
E'_2R_{21}E_1=R_{21}E'_2E_1\,.
$$

To complete the proof we have to consider, on the contrary, the algebra
generated by the entries of blocks $\Zm$, $1\le m\le N-1$, and determined by
the relations (3.22), (3.23) and to interpret the equality (3.21) as the
defining relation for \sZ. Clearly, \sZ defined this way is upper triangular
with units on the diagonal and
$$
\Em\sZa=\sum_{n=1}^m \Fnm(\En+\gZn)=\sZa(\Em+\gZm)\,,
$$
by (3.22) and (3.25) ($\Fnm\Em=0$ for $n\ge m+1$). Thus it remains to show
that this matrix \sZ verifies also (3.3). Using (3.17) one derives
\begin{eqnarray*}
R_{12}\sZa_1\diag(R)\sZa_2 & = &
\sum_{m=1}^N\sum_{n=1}^N R_{12}\Fmm_1(\Em_1+\gZm_1)\diag(R)
\Fnm_2(\En_2+\gZn_2)  \cr
& = & \sum_{m=1}^N\sum_{n=1}^N R_{12}\Fmm_1\Fnm_2(\Em_1+\gZm_1)R_{12}
(\En_2+\gZn_2)\,.  \cr
\end{eqnarray*}
Split the last double some into two pieces according to whether $m\le n$ or
$m>n$. Using (3.23) or its equivalent
$$
m\ge n \Longrightarrow
R_{12}^{-1}(\Em_1+\gZm_1)R_{12}(\En_2+\gZn_2)=(\En_2+\gZn_2)R_{21}
(\Em_1+\gZm_1)R_{21}^{-1}
$$
one finds that this double sum equals
\begin{eqnarray*}
& & \sum_{\quad m\le n}\!\!\!\!\!\sum R_{12}\Fmm_1\Fnm_2R_{21}^{-1}
(\En_2+\gZn_2)R_{21}(\Em_1+\gZm_1)R_{12} \cr
& & \qquad +\sum_{\quad m>n}\!\!\!\!\!\sum
R_{12}\Fmm_1\Fnm_2R_{12}(\En_2+\gZn_2)
R_{21}(\Em_1+\gZm_1)R_{21}^{-1} \cr
& & \ \ = \sum_m \Fmm_2(\Em_2+\gZm_2)\diag(R)\Fmm_1(\Em_1+\gZm_1)R_{12} \cr
& & \qquad +\sum_{\quad m<n}\!\!\!\!\!\sum (R_{21}^{-1}+(q-q^{-1})P)
\Fnm_2(\En_2+\gZn_2)\diag(R)\Fmm_1(\Em_1+\gZm_1)R_{12} \cr
& & \qquad +\sum_{\quad m>n}\!\!\!\!\!\sum
\Fnm_2(\En_2+\gZn_2)\diag(R)\Fmm_1(\Em_1+\gZm_1)
(R_{12}-(q-q^{-1})P) \cr
& & \ \ = \sZa_2\diag(R)\sZa_1R_{12}\,. \cr
\end{eqnarray*}
Here we have used repeatedly (3.17), (3.18), (3.19); for example,
\begin{eqnarray*}
m<n & \Longrightarrow &
P\Fnm_2(\En_2+\gZn_2)\diag(R)\Fmm_1(\Em_1+\gZm_1)R_{12} \cr
& & = \Fnm_1\Fmm_2(\En_1+\gZn_1)R_{12}(\Em_2+\gZm_2)R_{21}P \cr
& & =  \Fnm_1\Fmm_2(\Em_2+\gZm_2)\diag(R)(\En_1+\gZn_1)P \,,
\end{eqnarray*}
and
$$
m<n \Longrightarrow R_{21}^{-1}\Fmm_1\Em_1\Fnm_2\En_2 =
\Fmm_1\Em_1\Fnm_2\En_2\,.\quad\QED
$$

\section{A CONSTRUCTION OF LEFT \lUh -- \newline
MODULES}

The goal of this section is to equip the algebra \Fah with the structure
of a left \Uh -- module depending on $(N-1)$ parameters. More precisely,
here we will use the FRT description and hence deal with
the algebra \Ad (more or less equivalent to \AqAN) introduced in Section 2
(cf. (2.10)). The module will be redefined in terms of Chevalley generators
in Section 6. Moreover, the number of parameters is $N$ for the
\Ad--module but the expressions in Chevalley generators will turn out
to depend only on their ratios and so the number is reduced by one (as it
should be).

Assume for a moment that we are able to quantize the big cell of the flag
manifold with its real analytic structure, that means to construct
a $\ast$-algebra \cF generated jointly by $\zeta_{jk}$ and
$\zeta^\ast_{st}$ ($1\le j<k\le N,\ 1\le s<t\le N$) and containing both
\Fh and \Fah as its subalgebras. The strategy would be then to find
a (restriction) morphism $\psi:\Ada\to{\cal F}$ giving \cF the structure
of a left \Ad--module and afterwards to factorize off the ``holomorphic''
generators $\zeta_{jk}$ and to identify the factor-space
${\cal F}/\langle\zeta_{jk}\rangle$ with \Fah. The morphism $\psi$ is
determined by its values on the generators and thus it suffices to prescribe
a matrix $\psi(M)\in\Mat(N,{\cal F})$ obeying the corresponding relations
(2.10). In what follows we shall abuse the notation by writing simply $M$
instead of $\psi(M)$. However, as observed already in Ref. \cite{JALG},
it is not necessary to know the structure of the algebra \cF in full detail
for successful derivation of the left module structure. A left \Ad--module
on \Fah will be defined by prescribing the action on the unit and then
extending it to all polynomials in non-commutative variables $\zeta^\ast_{st}$
with the help of a recursive rule.

Ignoring the fact that we don't know the commutation relations between
the holomorphic and antiholomorphic generators we will introduce
projector-like matrices, one for each Grassmann manifold ($1\le m\le N-1$),
$$
\gQm:=\left(\begin{array}{c} \bI \\ \Zmast \end{array}\right)
(\bI+\Zm\Zmast)^{-1}
\left(\begin{array}{cc} \bI & \Zm \end{array}\right)\,.
\eqno(4.1)$$
Note that
$$
\bI-\gQm=\left(\begin{array}{c} -\Zm \\ \bI \end{array}\right)
(\bI+\Zmast\Zm)^{-1}
\left(\begin{array}{cc} -\Zmast & \bI \end{array}\right)\,.
\eqno(4.2)$$
We set also whenever convenient $\gQ^{(0)}:=\bO,\ \gQ^{(N)}:=\bI$. It holds
true that
$$
{\gQm}^\ast=\gQm\quad and\quad \gQm\gQn=\gQn\gQm=\gQ^{(\min(m,n))}\,.
\eqno(4.3)$$
The second equality means that $(\bI-\gQm)\gQn=0$ for $m\le n$ and follows
immediately from (4.1), (4.2) and (3.24). A morphism $\psi:\Ada\to\cFa$ is
prescribed by postulating a diagonalization of the quantum positive matrix
$M=\Lambda^\ast\Lambda$,
$$
\psi(M)=\xi_1\gQ^{(1)}+\xi_2(\gQ^{(2)}-\gQ^{(1)})+\dots+
\xi_N(\bI-\gQ^{(N)})\,.
$$
The construction of a left \Ad--module is based on the requirement that this
prescription actually defines a homomorphism.
\proclaim Lemma 4.1.
Assume that $N$ real (or real like, i.e., from $\Ra[[h]]\subset\Ca[[h]]$)
and mutually different parameters $\xi_m,\ 1\le m\le N$, are given.
The matrix
$$
M=\sum_{m=1}^N \xi_m(\gQm-\gQmm)
\eqno(4.4)$$
verifies the relations (repeating (2.10))
$$
M_2R_{12}^{-1}M_1R_{21}^{-1}=R_{12}^{-1}M_1R_{21}^{-1}M_2,\quad
M^\ast=M\,,
\eqno(4.5)$$
if and only if it holds
$$
\gQm_2R_{12}^{-1}M_1R_{21}^{-1}=R_{12}^{-1}M_1R_{21}^{-1}\gQm_2,\quad
1\le m\le N-1\,.
\eqno(4.6)$$

\noindent{\em Proof}.
Clearly $M^\ast=M$. Note that (4.3) means
$$
(\gQm-\gQmm)(\gQn-\gQnm)=\delta_{mn}(\gQm-\gQmm)\,.
$$
Multiply the first equality in (4.5) by $(\gQm-\gQmm)$ from the left and by
$(\gQn-\gQnm)$ from the right to obtain
$$
(\xi_m-\xi_n)(\gQm-\gQmm)R_{12}^{-1}M_1R_{21}^{-1}(\gQn-\gQnm)=0\,.
$$
The factor $(\xi_m-\xi_n)$ can be canceled for $m\not=n$. The summation
over $n$, with $m$ being fixed, results in
$$
(\gQm-\gQmm)R_{12}^{-1}M_1R_{21}^{-1}=
(\gQm-\gQmm)R_{12}^{-1}M_1R_{21}^{-1}(\gQm-\gQmm)\,.
$$
Here the RHS is self-adjoint ($R_{12}^\ast=R_{21}$)
and so the same is true for the LHS. This means that
$$
(\gQn-\gQnm)R_{12}^{-1}M_1R_{21}^{-1}=
R_{12}^{-1}M_1R_{21}^{-1}(\gQn-\gQnm)\,.
\eqno(4.7)$$
To get (4.6) it suffices to sum the equalities (4.7) over $n$, $1\le n\le m$.
The opposite implication is obvious since (4.6) means that
$\gQm_2,\ 0\le m\le N$, commute with $R_{12}^{-1}M_1R_{21}^{-1}$ and
so does $M_2$.
\QED

The commutation relation between $\gQm$ and $M$ implies a commutation
relation between $\gZmast$ and $M$. This will provide us with the desired
recursive rule (stated in (4.27) below).
\proclaim Lemma 4.2.
The commutation relation (4.6) is equivalent to
$$
(\bI-\Fm_2R_{12}\gZmast_2R_{12}^{-1})M_1R_{21}^{-1}(\Em_2+\gZmast_2)=
\Em_2M_1R_{21}^{-1}\,.
\eqno(4.8)$$
Furthermore, the matrix $(\bI-\Fm_2R_{12}\gZmast_2R_{12}^{-1})$ is
invertible and it holds true that
\begin{eqnarray*}
(\bI-\Fm_2R_{12}\gZmast_2R_{12}^{-1})^{-1}\Em_2
& = & R_{12}(\Em_2+\gZmast_2)R_{21}
\makebox[1em]{}\cr
&   & -(q-q^{-1})(\Em_1+\gZmast_1)R_{12}P(\Em_1+\gZmast_1)\,.
\makebox[0em][r]{}\cr
&   &
\!\makebox[21em][r]{(4.9)}\cr
\end{eqnarray*}

\noindent{\em Proof}.
In this proof we shall suppress the superscript $(m)$.
(a) The relation (4.6) means that
$$
(\bI-\gQ_2)R_{12}^{-1}M_1R_{21}^{-1}\gQ_2=0\,.
$$
From the structure of $\gQ$ (cf. (4.1), (4.2)) one finds that this is
equivalent to
$$
(\bI-E_2-\gZast_2)R_{12}^{-1}M_1R_{21}^{-1}(E_2+\gZast_2)=0
$$
and hence
$$
(F_2-\gZast_2)R_{12}^{-1}M_1R_{21}^{-1}\gZast_2 =
-(F_2-\gZast_2)R_{12}^{-1}M_1R_{21}^{-1}E_2\,.
\eqno(4.10)$$
Recall that $\gZast=F\gZast E$. Thus the both sides of (4.10) are invariant
with respect to multiplication by $F_2$ from the left. Since (cf. (2.15))
$$
F_2R_{12}^{\pm1}F_2R_{12}^{\mp1}=F_2\,,
$$
the multiplication by $F_2R_{12}$ from the left results in an equivalent
equality, namely
$$
(F_2-F_2R_{12}\gZast_2R_{12}^{-1})M_1R_{21}^{-1}\gZast_2 =
-(F_2-F_2R_{12}\gZast_2R_{12}^{-1})M_1R_{21}^{-1}E_2\,.
\eqno(4.11)$$
Now it suffices to observe that $F_2R_{21}^{-1}\gZast_2=R_{21}^{-1}\gZast_2$
and $F_2R_{21}^{-1}E_2=R_{21}^{-1}E_2-E_2R_{21}^{-1}$ in order to conclude
that (4.11) coincides with (4.8).

(b) In the lexicographically ordered basis, $R_{12}^{\pm1}$ is lower
triangular ($R_{jk,st}=0$ for $j<s$ and
$R_{jk,jt}=\delta_{kt}R_{jk,jt}$). It follows readily that
$F_2R_{12}\gZast_2R_{12}^{-1}$ is lower triangular with vanishing diagonal
and hence nilpotent and consequently $(\bI-F_2R_{12}\gZast_2R_{12}^{-1})$
is invertible. It remains to show that
$$
E_2 = (\bI-F_2R_{12}\gZast_2R_{12}^{-1})[R_{12}(E_2+\gZast_2)R_{21}
-(q-q^{-1})(E_1+\gZast_1)R_{12}P(E_1+\gZast_1)]\,.
$$
We have for the LHS (cf. (2.15), (2.4), (2.17))
\begin{eqnarray*}
E_2 & = & E_2R_{21}^{-1}(E_2+\gZast_2)R_{21} \cr
& = & E_2R_{12}(E_2+\gZast_2)R_{21}
-(q-q^{-1})E_2(E_1+\gZast_1)R_{12}P(E_1+\gZast_1)
\end{eqnarray*}
while the RHS can be expanded and treated using (3.23) and
$$
F_2R_{12}\gZast_2=F_2R_{12}(E_2+\gZast_2),\quad
(E+\gZast)^2=(E+\gZast)\,.
$$
This gives immediately the result.
\QED

In addition to the recursive rule we also need to know the action on the
unit, i.e., the expression $M\cdot1$. Thinking of \Fah as the factor space
${\cal F}/\langle\zeta_{jk}\rangle$ we require naturally that
$\Zm\cdot1=0,\ \forall m$. This implies
$$
\gQm\cdot1=\left(\begin{array}{c} \bI \\ \Zmast \end{array}\right)
(\bI+\Zm\Zmast)^{-1}
\left(\begin{array}{cc} \bI & \bO \end{array}\right) \cdot1\,.
$$
Guided by the experience obtained when working with quantum Grassmannians
\cite{QG} we accept as an ansatz that for some scalar $\eta_m$,
$$
(\bI+\Zm\Zmast)^{-1}\cdot1=\eta_m\bI\,,
$$
whence
$$
\gQm\cdot1=\eta_m(\Em+\gZmast)\,.
\eqno(4.12)$$
To get rid of superfluous parameters we use the substitution
$$
\lambda_N:=\xi_N,\
\lambda_{N-j}:=\lambda_{N-j+1}+(\xi_{N-j}-\xi_{N-j+1})\eta_{N-j}\quad
for\ j=1,\dots,N-1\,.
\eqno(4.13)$$

Now we are ready to state the result provided \Fah is described in terms of
\Gmah, $1\le m\le N-1$, as given in Proposition 3.3. Before let us
formulate some auxiliary relations needed for the proof. At the same time,
we introduce some shorthand notation, also only for the purpose of this
proof.  Set
\begin{eqnarray*}
\makebox[8em]{}
\gamma & := & q-q^{-1}\,,
\!\makebox[20em][r]{(4.14)}\cr
\makebox[8em]{}
\gXm & := & \Em+\gZmast\,,
\,\,\makebox[17em][r]{(4.15)}\cr
\makebox[8em]{}
\Xm_{12} & := & (\bI-\Fm_2R_{12}\gZmast_2R_{12}^{-1})^{-1}\Em_2\,,
\makebox[10em][r]{(4.16)}\cr
\end{eqnarray*}
and $X^{(0)}_{12}:=\bO,\ X^{(N)}_{12}:=\bI$. Thus, by Proposition 3.3 and
(3.25), we have
\begin{eqnarray*}
\makebox[7em]{}
m\le n & \Longrightarrow & \gXn\gXm=\gXm\gXn=\gXm\,,
\!\!\makebox[11em][r]{(4.17)}\cr
\makebox[7em]{}
m\le n & \Longrightarrow & R_{21}\gXm_1R_{12}\gXn_2=\gXn_2R_{21}\gXm_1R_{12}\,,
\!\makebox[8em][r]{(4.18)}\cr
\end{eqnarray*}
and by Lemma 4.2,
$$
\Xm_{12}=R_{12}\gXm_2R_{21}-\gamma\gXm_1R_{12}P\gXm_1\,.
\eqno(4.19)$$
\proclaim Lemma 4.3.
For $1\le m\le n\le N$, it holds true that
\begin{eqnarray*}
\makebox[5em]{}
\Xn_{12}\Xm_{12} & = & \Xm_{12}\,,
\,\makebox[14em][r]{(4.20)}\cr
\makebox[5em]{}
R_{32}\Xm_{12}R_{23}\Xn_{13} & = & \Xn_{13}R_{32}\Xm_{12}R_{23}\,,
\,\makebox[9em][r]{(4.21)}\cr
\makebox[5em]{}
(\Xm_{21}\gXn_2+\Xn_{21}\gXm_2)R_{12}^{-1} & = &
R_{12}^{-1}(\Xm_{12}\gXn_1+\Xn_{12}\gXm_1)\,.
\makebox[5em][r]{(4.22)}\cr
\end{eqnarray*}

\noindent{\em Proof}.
(4.20): Observe that
$$
\Xm_{12}=R_{12}\gXm_2R_{21}(\bI-\gamma\gXm_1R_{21}^{-1}P_{12})
\eqno(4.23)$$
and ($m\le n$)
$$
(\bI-\gamma\gXn_1R_{21}^{-1}P_{12})R_{12}\gXm_2R_{21}=R_{21}^{-1}
\gXm_2R_{21}\,.
$$
Consequently,
$$
\Xn_{12}\Xm_{12} = R_{12}\gXn_2R_{21}\cdot R_{21}^{-1}\gXm_2R_{21}
(\bI-\gamma\gXm_1R_{21}^{-1}P_{12})  = \Xm_{12}\,.
$$

(4.21): Using repeatedly YB equation, (4.18) and (4.17), one can
derive that
\begin{eqnarray*}
\gXn_3R_{13}^{-1}R_{32}(R_{12}\gXm_2R_{21})R_{23}R_{13}\gXn_3 & = &
\gXn_3R_{12}R_{32}\gXm_2R_{23}R_{21}\gXn_3 \\
& = & R_{13}^{-1}R_{32}(R_{12}\gXm_2R_{21})R_{23}R_{13}\gXn_3\,,
\end{eqnarray*}
and (making use also of $R_{13}^{-1}=R_{31}-\gamma P_{13}$)
\begin{eqnarray*}
\gXn_3R_{13}^{-1}R_{32}(\gXm_1R_{12}P_{12}\gXm_1)R_{23}R_{13}\gXn_3 &=&
\gXn_3R_{13}^{-1}\gXm_1R_{13}P_{12}R_{21}R_{31}\gXm_1R_{13}\gXn_3 \\
&=& R_{13}^{-1}R_{32}(\gXm_1R_{12}P_{12}\gXm_1)R_{23}R_{13}\gXn_3\,.
\end{eqnarray*}
In virtue of (4.19), this means jointly that
$$
\gXn_3R_{13}^{-1}R_{32}\Xm_{12}R_{23}R_{13}\gXn_3=
R_{13}^{-1}R_{32}\Xm_{12}R_{23}R_{13}\gXn_3\,,
$$
and consequently, utilizing once more (4.23),
$$
(\Fn_3-\Fn_3R_{13}\gXn_3R_{13}^{-1})R_{32}\Xm_{12}R_{23}\Xn_{13}=0\,.
\eqno(4.24)$$
Furthermore,
\begin{eqnarray*}
\makebox[7em]{}
\En_3R_{32}\Xm_{12}R_{23}\Xn_{13} & = &
\En_3R_{32}\Xm_{12}\En_2\En_3R_{23}\Xn_{13}
\makebox[6em][r]{}\cr
\makebox[7em]{}
& = & \En_3R_{32}\Xm_{12}R_{23}\,,
\makebox[10em][r]{(4.25)}\cr
\end{eqnarray*}
since $\En_3\Xn_{13}=\En_3$. Summing (4.24) and (4.25) we get the sought
equality in the form ($\Fn_3R_{13}\gXn_3=\Fn_3R_{13}\gZnast_3$)
$$
(\bI-\Fn_3R_{13}\gZnast_3R_{13}^{-1})R_{32}\Xm_{12}R_{23}\Xn_{13}=
\En_3R_{32}\Xm_{12}R_{23}\,.
$$

(4.22): Using again (4.19), (4.17) and (4.18) as well as (2.4) we find that
\begin{eqnarray*}
\LHS & = & R_{21}\gXm_1R_{12}\gXn_2R_{12}^{-1}
-\gamma\gXm_2R_{21}P\gXm_2R_{12}^{-1} \cr
& & \ +R_{21}\gXn_1R_{12}\gXm_2R_{12}^{-1}
-\gamma\gXn_2R_{21}P\gXm_2R_{12}^{-1} \cr
& = & \gXn_2R_{21}\gXm_1-\gamma R_{12}^{-1}\gXm_1R_{12}\gXm_2P
+\gXm_2R_{21}\gXn_1R_{21}^{-1}R_{12}^{-1}\,, \cr
\end{eqnarray*}
and
$$
\RHS = \gXm_2R_{21}\gXn_1-\gamma R_{12}^{-1}\gXm_1R_{12}P\gXm_1
+\gXn_2R_{21}\gXm_1-\gamma R_{12}^{-1}\gXn_1R_{12}P\gXm_1\,.
$$
Consequently,
$$
\LHS-\RHS =\gXm_2R_{21}\gXn_1R_{21}^{-1}(R_{12}^{-1}-R_{21})
+ \gamma R_{12}^{-1}\gXn_1R_{12}\gXm_2P  = 0\,.
\quad\QED
$$
\proclaim Proposition 4.4.
The relations
\begin{eqnarray*}
&  & \!\!\!\!\!\!\!\!
M\cdot1=(\lambda_1-\lambda_2)(E^{(1)}+{\gZ^{(1)}}^\ast)+\dots+
(\lambda_{N-1}-\lambda_N)(E^{(N-1)}+{\gZ^{(N-1)}}^\ast)+\lambda_N\bI\,, \cr
&& \,\makebox[34em][r]{(4.26)}\cr
&   & \!\!\!\!\!\!\!\!
M_1R_{21}^{-1}\cdot(\Em_2+\gZmast_2)f=
(\bI-\Fm_2R_{12}\gZmast_2R_{12}^{-1})^{-1}\Em_2M_1R_{21}^{-1}\cdot f\,,
\,\,\makebox[5em][r]{(4.27)}\cr
\end{eqnarray*}
for all $f\in\mbox{\Fah}$, $1\le m\le N-1$, define on \Fah unambiguously the
structure of a left \Ad--module (the central dot ``$\cdot$'' stands for the
action) depending on $N$ scalar parameters $\lambda_1,\dots,\lambda_N$.

\noindent{\em Proof}.
As already mentioned, the idea of defining a module this way was utilized
in Ref. \cite{JALG} (Proposition 5.4) and the proof is quite similar, too.
Here we rely heavily on Proposition 3.3 giving a description of \Fah in
terms of \Gmah.
First we have to show that (4.26) and (4.27) define correctly a linear mapping
$\Faha\to\Mat(N,\Faha) : f\mapsto M\cdot f$. Let $\bF$ be the free
algebra generated by $\bzm_{jk},\ 1\le j\le m<k\le N$,
and the matrix $\bgZm$ be obtained from $\gZmast$ when replacing the entries
$\zmast_{jk}$ by $\bzm_{jk}$. Furthermore, we use in an obvious sense the
symbols $\bgXm$ and $\bXm_{12}$ parallelly to (4.15) and (4.16), respectively.
Hence \Fah is obtained from $\bF$
by means of factorization by the two-sided ideal generated by the
relations (4.17), (4.18), with $\gXm$'s being replaced by $\bgXm$'s,
and the elements $\zmast_{jk}$ are the factor images of
$\bzm_{jk}$.
Doing the same replacement in (4.26) and (4.27)
one obtains a well defined linear mapping
$$
\bF\to\Mat(N,\bF): f\mapsto M\cdot f .
\eqno(4.28) $$
A straightforward calculation gives for $m\le n$ and $\forall f\in\bF$,
\begin{eqnarray*}
& & M_1R_{21}^{-1}\cdot(\bgXn_2\bgXm_2-\bgXm_2)f=
(\bXn_{12}\bXm_{12}-\bXm_{12})M_1R_{21}^{-1}\cdot f\,, \cr
& & M_1R_{31}^{-1}R_{21}^{-1}\cdot
(R_{32}\bgXm_2R_{23}\bgXn_3-\bgXn_3R_{32}\bgXm_2R_{23})f \cr
& & \quad = (R_{32}\bXm_{12}R_{23}\bXn_{13}-\bXn_{13}R_{32}\bXm_{12}R_{23})
M_1R_{31}^{-1}R_{21}^{-1}\cdot f\,. \cr
\end{eqnarray*}
This means that the linear mapping (4.28) factorizes from $\bF$ to
\Fah if and only if the factor-images of the matrices
$(\bXn_{12}\bXm_{12}-\bXm_{12})$ and
$(R_{32}\bXm_{12}R_{23}\bXn_{13}-\bXn_{13}R_{32}\bXm_{12}R_{23})$
vanish. But these are exactly the relations (4.20) and (4.21) proven in
Lemma 4.3.

To show that \Fah is really a left \Ad--module we have to verify
the equality
$$
M_2R_{12}^{-1}M_1R_{21}^{-1}\cdot1=R_{12}^{-1}M_1R_{21}^{-1}M_2\cdot1
\eqno(4.29) $$
and the implication
\begin{eqnarray*}
& & (M_2R_{12}^{-1}M_1R_{21}^{-1})R_{31}^{-1}R_{32}^{-1}\cdot f=
(R_{12}^{-1}M_1R_{21}^{-1}M_2)R_{31}^{-1}R_{32}^{-1}\cdot f\Longrightarrow
\cr
& & (M_2R_{12}^{-1}M_1R_{21}^{-1})R_{31}^{-1}R_{32}^{-1}\cdot\gZmast_3f=
(R_{12}^{-1}M_1R_{21}^{-1}M_2)R_{31}^{-1}R_{32}^{-1}\cdot\gZmast_3f\,,
\makebox[4em][r]{(4.30)}\cr
\end{eqnarray*}
for $\forall f\in\Faha,\ 1\le m\le N-1$,
since then (4.29) and (4.30) jointly imply
$$
M_2R_{12}^{-1}M_1R_{21}^{-1}\cdot f=R_{12}^{-1}M_1R_{21}^{-1}M_2\cdot f,
\quad\forall f\in\Faha\, .
$$

{\it Verification of (4.29)}.
Using (4.26) and (4.27) one finds easily that (4.29) means ($\lambda_{N+1}:=0$)
\begin{eqnarray*}
& & \sum_{m=1}^N(\lambda_m-\lambda_{m+1})\Xm_{21}
\sum_{n=1}^N(\lambda_n-\lambda_{n+1})\gXn_2R_{12}^{-1}R_{21}^{-1} \cr
& & = R_{12}^{-1}\sum_{m=1}^N(\lambda_m-\lambda_{m+1})\Xm_{12}
\sum_{n=1}^N(\lambda_n-\lambda_{n+1})\gXn_1R_{21}^{-1}\,. \cr
\end{eqnarray*}
But this equality follows immediately from the relation (4.22) proven in
Lemma 4.3.

{\it Verification of (4.30)}.
Note first that the equality after the sign of implication in (4.30) can be
replaced by
$$
(M_2R_{12}^{-1}M_1R_{21}^{-1})R_{31}^{-1}R_{32}^{-1}\cdot\gXm_3f=
(R_{12}^{-1}M_1R_{21}^{-1}M_2)R_{31}^{-1}R_{32}^{-1}\cdot\gXm_3f\,,
$$
Furthermore, reversing the proof of Lemma 4.2, part (a), one finds that (4.27)
is equivalent to
$$
R_{12}^{-1}M_1R_{21}^{-1}\cdot\gXm_2f=
\gXm_2R_{12}^{-1}M_1R_{21}^{-1}\cdot\gXm_2f,\quad\forall f\in\Faha\,.
\eqno(4.31)$$
Using repeatedly YB equation and (4.31) one derives that
\begin{eqnarray*}
\makebox[4em]{}
& & (\bI-\gXm_3)R_{23}^{-1}R_{13}^{-1}(M_2R_{12}^{-1}M_1R_{21}^{-1})
R_{31}^{-1}R_{32}^{-1}\cdot\gXm_3f \cr
& & = (\bI-\gXm_3)(R_{23}^{-1}M_2R_{32}^{-1})R_{12}^{-1}
(R_{13}^{-1}M_1R_{31}^{-1})\cdot\gXm_3R_{21}^{-1}f \cr
& & =0\,,
\makebox[28em][r]{(4.32)}\cr
\end{eqnarray*}
for $(\bI-\gXm_3)\gXm_3=0$. It follows from (4.32) that
\begin{eqnarray*}
\makebox[4em]{}
& & (\bI-\Fm_3R_{13}R_{23}\gXm_3R_{23}^{-1}R_{13}^{-1})
(M_2R_{12}^{-1}M_1R_{21}^{-1})R_{31}^{-1}R_{32}^{-1}\cdot\gXm_3f \cr
& & = \Em_3(M_2R_{12}^{-1}M_1R_{21}^{-1})R_{31}^{-1}R_{32}^{-1}\cdot f\,,
\makebox[15em][r]{(4.33)}\cr
\end{eqnarray*}
for $\Em_3R_{31}^{-1}R_{32}^{-1}\gXm_3=\Em_3R_{31}^{-1}R_{32}^{-1}$.
Quite similarly one obtains
\begin{eqnarray*}
\makebox[4em]{}
& & (\bI-\Fm_3R_{13}R_{23}\gXm_3R_{23}^{-1}R_{13}^{-1})
(R_{12}^{-1}M_1R_{21}^{-1}M_2)R_{31}^{-1}R_{32}^{-1}\cdot\gXm_3f \cr
& & = \Em_3(R_{12}^{-1}M_1R_{21}^{-1}M_2)R_{31}^{-1}R_{32}^{-1}\cdot f\,,
\makebox[15em][r]{(4.34)}\cr
\end{eqnarray*}
The right hand sides of (4.33) and (4.34) are equal by assumption and so,
to complete the proof, it suffices to show that that the matrix
$(\bI-\Fm_3R_{13}R_{23}\gXm_3R_{23}^{-1}R_{13}^{-1})$
is invertible. But using a similar argument as in the proof of Lemma 4.2,
part (b), one finds that the matrix
$$
\Fm_3R_{13}R_{23}\gXm_3R_{23}^{-1}R_{13}^{-1}=
\Fm_3R_{13}R_{23}\gZmast_3R_{23}^{-1}R_{13}^{-1}
$$
is nilpotent. \QED

\section{LEIBNIZ RULE AND THE DRESSING\newline
TRANSFORMATION}

The right dressing transformation
$$
\cR:\Aqa(AN)\to\Aqa(AN)\otimes\Aqa(SU(N))
$$
can be introduced formally using the canonical element in
$\Aqa(AN)\otimes\Aqa(SU(N))$ \cite{J-S-I, J-S-II} and it factorizes from
\AqAN to both \Fh and \Fah. Dually it induces a left action of \Uh on
\AqAN (or \Fh or \Fah) via the pairing $\langle\cdot,\cdot\rangle$ between
\Uh and \AqSU,
$$
\xi_Y\cdot f:=(\id\otimes\langle Y,\cdot\rangle)\cR(f)\,.
\eqno(5.1)$$
Leibniz rule for $\xi$ means that
$$
\xi_Y\cdot fg=(\xi_{Y_{(1)}}\cdot f)(\xi_{Y_{(2)}}\cdot g),\quad
with\ \DeltaYa\,.
\eqno(5.2)$$
The aim of this section is to show that Leibniz rule for $\xi$ (acting on
\Fah) induces the recursive rule (4.27).

Here we introduce the dressing transformation
$$
\cR:\Fha\to\Fha\otimes\Aqa(SU(N))
\eqno(5.3)$$
directly by prescribing its values on the generators $\zeta_{jk}$ arranged
in the matrix \sZ. As usual, the identifications $\Fha\equiv\Fha\otimes1$
and $\Aqa(SU(N))\equiv1\otimes\Aqa(SU(N))$ simplify a lot the notation.
We define, with the help of Gauss decomposition and formally in the same
way as classically,
$$
\cR(\sZa):=(\sZa U)_+,\quad where\ \sZa U =(\sZa U)_-(\sZa U)_+
\eqno(5.4)$$
and $(\sZa U)_+$ is upper triangular with units on the diagonal
while $(\sZa U)_-$ is lower triangular ($U$ still designates the
vector corepresentation of $\Aqa(SU(N))$\ ). $\cR$ extends to an algebra
homomorphism and it holds
\begin{eqnarray*}
\makebox[11em]{}
(\id\otimes\vare)\circ\cR & = & \id\,,
\makebox[16em][r]{(5.5)}\cr
\makebox[11em]{}
(\cR\otimes\id)\circ\cR & = & (\id\otimes\Delta)\circ\cR\,,
\,\,\makebox[11em][r]{(5.6)}\cr
\end{eqnarray*}
as follows readily from (2.1) and the uniqueness of Gauss decomposition.

Let us rewrite the dressing transformation (5.4) in terms of coordinate
functions $\zm_{jk}$ on quantum Grassmannians. From (3.20) on finds that
\begin{eqnarray*}
\Em\sZa U\cR(\Em+\gZm) & = & \Em\sZa U\cR(\sZa)^{-1}\Em\cR(\sZa) \cr
& = & \Em(\sZa U)_- (\bI-\Fm)(\sZa U)_+ \cr
& = & \Em\sZa U\,, \cr
\end{eqnarray*}
for $\Em(\sZa U)_-\Fm=0$ owing to the fact that  $(\sZa U)_-$ is lower
triangular. Similarly, $\Fm\cR(\sZa)^{-1}\Em=0$ and thus
$$
(\sZa\Fm+\Em\sZa U)\cR(\Em+\gZm)=\Em\sZa U\,.
$$
Consequently,
\begin{eqnarray*}
\makebox[3em]{}
\cR(\Em+\gZm) & = & (\Fm+\sZa^{-1}\Em\sZa U)^{-1}\sZa^{-1}\Em\sZa U
\makebox[4em][r]{}\cr
\makebox[3em]{}
& = & [\Fm+(\Em+\gZm)U]^{-1}(\Em+\gZm)U\,.
\makebox[5em][r]{(5.7)}\cr
\end{eqnarray*}
This can be expressed also in terms of blocks $\Zm$. Decompose $U$ into
blocks,
$$
U=\left(\begin{array}{cc} \Am & \Bm \\ \Cm & \Dm \end{array}\right)\,,
$$
where $\Am$ has dimension $m\times m$, $\Bm$ has dimension $m\times(N-m)$
etc. The result is
$$
\cR(\Zm)=(\Am+\Zm\Cm)^{-1}(\Bm+\Zm\Dm)
$$
and this is the correct formula \cite{QG}. The dressing transformation on
\Fah is obtained readily by taking the adjoints. Particularly,
$$
\cR(\Em+\gZmast)=U^\ast(\Em+\gZmast)[\Fm+U^\ast(\Em+\gZmast)]^{-1}\,.
\eqno(5.8)$$

Let us add a couple of remarks. Despite of the fact that the formulas (5.4),
(5.7) contain rational singularities (a consequence of the localization from
the dressing orbit to the big cell) the action $\xi$ (local in its nature)
is free of any singularities. Furthermore, observe from (4.1) and
$$
\cR\bigl(\ \bI \quad \Zm\ \bigr) =
\bigl(\Am+\Zm\Cm\bigr)^{-1}
\bigl(\ \bI \quad \Zm\ \bigr)U
$$
that $\cR(\gQm)=U^\ast\gQm U$ and thus the quantum diagonalization
(4.4) is in accordance with the rule
$\cR(\Lambda^\ast\Lambda)=U^\ast\Lambda^\ast\Lambda U$, as it should be
\cite{J-S-I, QG}. Finally, recall that the comultiplication $\Delta$ is
always assumed to come from \Uh rather than from \AqAN. Now we can
formulate a statement relating Leibniz rule to the \Ad--module \Fah which
was defined in Proposition 4.4.
\proclaim Proposition 5.1.
For $\forall Y\in\Ada$ and $\forall \psi, f\in\Faha$, it holds
$$
Y\cdot \psi f=(\xi_{Y_{(1)}}\cdot\psi)(Y_{(2)}\cdot f),\quad
with\ \DeltaYa\,.
\eqno(5.9)$$

{\em Remark}. This is in accordance with the classical case. The method of
orbits, if expressed in local coordinates, yields a representation of Lie
algebra {\gt su}$(N)$ in terms of first order differential operators,
$$
Y\in\mbox{\gt su}(N)\mapsto p_Y(\bar z)\partial_{\bar{z}}+q_Y(\bar z)\,,
$$
with generally non-vanishing zero-order term $q_Y(\bar z)$ and with the
first-order term identical to the vector field
$\xi_Y=p_Y(\bar z)\partial_{\bar{z}}$ coming from the infinitesimal
coadjoint action. Clearly,
$$
Y\cdot1=q_Y(\bar z)\quad and\quad
Y\cdot\psi f=(\xi_Y\cdot\psi)f+\psi(Y\cdot f)\,.
$$

{\em Proof}. As we are facing two actions ($\xi_X\xi_Y=\xi_{XY}$ and
$X\cdot(Y\cdot f)=(XY)\cdot f$) it is sufficient to verify (5.9) only for the
generators of \Ad, i.e., for $Y$ running over the entries of $M$. Note
that (5.9) holds trivially for $Y=1$ ($\Delta1=1\otimes1$) and the same is
true for $\psi=1$ owing to the equalities $\vare(Y_{(1)})Y_{(2)}=Y$ and
$$
\xi_Y\cdot1=(\id\otimes\langle Y,\cdot\rangle)1\otimes1=\vare(Y)\,.
$$
In virtue of Leibniz rule (5.2) and the co-associativity of $\Delta$ it
suffices, too, to verify (5.9) for $\psi$ running over the generators of \Fah.

Observe that
$$
M_2=\Lambda^\ast_2\Lambda_2=\tr_1(P_{12}\Lambda^\ast_1\Lambda_2)\,,
$$
where $\tr_1$ means the trace applied only in the first factor of the
tensor product in question. Thus we can formulate the problem in the
following way. Evaluate
$$
\Lambda^\ast_1\Lambda_2R_{32}^{-1}\cdot(\Em_3+\gZmast_3)f
$$
using the rule (5.9) and then apply to the obtained expression
$$
(\bI-\Fm_3R_{23}\gZnast_3R_{23}^{-1})\,\tr_1(P_{12}\,\cdot)\,.
$$
The result should coincide with that one given by the recursive rule (4.27),
i.e., with $\Em_3M_2R_{32}^{-1}\cdot f$. In the rest of the proof we drop
the superscript $(m)$.

As the comultiplication $\Delta_{AN}$ in (2.8) is opposite to
$\Delta$ we have
$$
\Delta\Lambda^t=\Lambda^t\dot{\otimes}\Lambda^t,\quad
\Delta\Lambda^\ast=\Lambda^\ast\dot{\otimes}\Lambda^\ast\,.
$$
Thus we start from (cf. (5.8))
\begin{eqnarray*}
\makebox[3em]{}
& & (\langle\Lambda^\ast_1\Lambda^t_2,\cdot\rangle \times
\Lambda^\ast_1\Lambda^t_2\cdot)^{t_2}
\left( R_{32}^{-1}U^\ast_3(E_3+\gZast_3)
[F_3+U^\ast_3(E_3+\gZast_3)]^{-1} \times f\right) \cr
\makebox[3em]{}
& & = \{(R_{32}^{-1})^{t_2}
\langle\Lambda^\ast_1\Lambda^t_2,\cdot\rangle U^\ast_3(E_3+\gZast_3)
[F_3+U^\ast_3(E_3+\gZast_3)]^{-1}\Lambda^\ast_1\Lambda^t_2\cdot f\}^{t_2}\,.
\makebox[3em][r]{(5.10)}\cr
\end{eqnarray*}
The pairing $\langle\Lambda^\ast_1\Lambda^t_2,\cdot\rangle$ acts on the
elements of algebra \AqSU occuring in the matrices $U^\ast_3$. In view of
(2.14), we have
$$
\langle\Lambda^\ast_1\Lambda^t_2,U^\ast_3\rangle=
\langle\Lambda^\ast_1\Lambda^t_2,U_3\rangle^{-1} =
((R_{32}^{-1})^{t_2})^{-1}R_{13}\,.
$$
Consequently, the expression (5.10) equals
$$
\{ R_{13}(E_3+\gZast_3)
[F_3+((R_{32}^{-1})^{t_2})^{-1}R_{13}(E_3+\gZast_3)]^{-1}
\Lambda^\ast_1\Lambda^t_2\cdot f\}^{t_2}\,.
$$
When applying $\tr_1(P_{12}\cdot)$ observe that (for any $X_{123}$)
$$
\tr_1\big( P_{12}(X_{123}\Lambda^\ast_1\Lambda^t_2)^{t_2}\big) =
\tr_1(P_{12}X_{123}^{t_2}M_1)
$$
and when multiplying by $(\bI-F_3R_{23}\gZast_3R_{23}^{-1})$ from the left
note that
\begin{eqnarray*}
(\bI-F_3R_{23}\gZast_3R_{23}^{-1})R_{23}(E_3+\gZast_3) & = &
R_{23}(E_3+\gZast_3)-F_3R_{23}(E_3+\gZast_3) \cr
& = & E_3R_{23}(E_3+\gZast_3)\,. \cr
\end{eqnarray*}
This way we arrive at the expression
$$
\tr_1\big( P_{12}\{E_3R_{13}(E_3+\gZast_3)
[F_3+((R_{32}^{-1})^{t_2})^{-1}R_{13}(E_3+\gZast_3)]^{-1}\}^{t_2}
M_1\cdot f\big)\,.
\eqno(5.11)$$
Finally we use (cf. (2.15))
$$
E_3=E_3(R_{32}^{-1})^{t_2}E_3((R_{32}^{-1})^{t_2})^{-1}
$$ and the obvious equality (multiply by $[\dots]$ from the right)
$$
E_3((R_{32}^{-1})^{t_2})^{-1}R_{13}(E_3+\gZast_3)
[F_3+((R_{32}^{-1})^{t_2})^{-1}R_{13}(E_3+\gZast_3)]^{-1}=E_3
$$
to conclude that (5.11) equals
\begin{eqnarray*}
\tr_1\big(P_{12}\{ E_3(R_{32}^{-1})^{t_2}E_3\}^{t_2}M_1\cdot f\big) & = &
E_3\,\tr_1(P_{12}R_{32}^{-1}M_1\cdot f) \cr
& = & E_3M_2R_{32}^{-1}\cdot f\,,
\end{eqnarray*}
as required. \QED

{\em Remark}. Let \tl be an irreducible finite-dimensional representation
of \Uh \linebreak
acting in a Hilbert space \Hl and corresponding to a lowest weight
$\lambda$ with a lowest weight vector \el.
Let $\Tla\in\End(\Hla)\otimes\Aqa(SU(N))$
be the related corepresentation of \AqSU. As shown in Ref.
\cite{J-S-II}, basically all information about the representation is encoded
in the element
$$
\wla:=\big(\,(\ela,\cdot\,\ela)\otimes\id\big)\Tla\in\Aqa(SU(N))\,.
$$
One can naturally realize the algebras \Fh and \Fah (but not the full
algebra $\cal F$ of real analytic functions) as subalgebras in the
localization of \AqSU when the element \wl is allowed to be invertible.
The comultiplication in \AqSU, being restricted to \Fah, coincides with
the dressing transformation. This means also that
$\Delta(\Faha)\subset\Faha\otimes\Aqa(SU(N))$. Dually we again get a left
action of \Uh on \Fah and keep the symbol $\xi$ for it. The formula
$$
(Y,f)\mapsto\wla^{-1}\xi_Y\cdot(\wla f)
$$
defines a left \Uh -- module structure on \Fah and the representation \tl
can be identified with the cyclic submodule \Ml with unit as the cyclic
vector and, at the same time, the lowest weight vector. With the same success
we could use the prescription
$$
(Y,f)\mapsto(\xi_Y\cdot(f\wla))\wla^{-1}\,.
\eqno(5.12)$$
It is not difficult to check that \tl can be again identified with a cyclic
submodule ${\cal M}'_\lambda$ with respect to this new action and the
unit is again the cyclic and the lowest weight vector. From (5.12) follows
easily that
\begin{eqnarray*}
Y\cdot \psi f & = & (\xi_Y\cdot(\psi f\wla))\wla^{-1} =
(\xi_{Y_{(1)}}\cdot\psi)(\xi_{Y_{(2)}}\cdot(f\wla))\wla^{-1} \cr
& = & (\xi_{Y_{(1)}}\cdot\psi)(Y_{(2)}\cdot f)\,, \cr
\end{eqnarray*}
with \DeltaY.

\section{IRREDUCIBLE REPRESENTATIONS}

The left action $\xi$, dual to the dressing transformation
$\cR:\Faha\to\Faha\otimes\Aqa(SU(N))$, can be expressed explicitly in terms
of Chevalley generators. For graphical reasons, we shall write
$\xi(Y)\cdot f$ instead of $\xi_Y\cdot f$. Relying on Proposition 5.1, one
can pass from the FRT description to Chevalley generators also in the
definition of the left module structure on \Fah. First note that the rules
(2.11), (2.7) imply
$$
\langle\qpmHj,A^{-1}\rangle=\langle\qpmHj,A\rangle^{-1},\
\langle\Xpmj,A^{-1}\rangle=-\langle q^{H_j/2},A\rangle^{-1}
\langle\Xpmj,A\rangle\langle q^{-H_j/2},A\rangle^{-1}\,,
\eqno(6.1)$$
where $A$ is any invertible square matrix with entries from \AqSU.
In particular (cf. (2.13)),
$$
\langle\qpmHj,U^\ast\rangle=q^{\mp(\Ejj-\Ejpjp)},\
\langle\Xpj,U^\ast\rangle=-q^{-1}\Ejjp,\
\langle\Xmj,U^\ast\rangle=-q\Ejpj\,.
\eqno(6.2)$$
To evaluate $\xi$ one can use the definition (5.1), the formula (5.8) for the
dressing transformation and the rules (2.11), (6.1), (6.2):
\begin{eqnarray*}
\makebox[3em]{}
\xi(\qpmHj)\cdot \gZmast & = &
\Fm q^{\mp(\Ejj-\Ejpjp)}(\Em+\gZmast) \cr
& & \times [\Fm+q^{\mp(\Ejj-\Ejpjp)}(\Em+\gZmast)]^{-1} \cr
& = & q^{\mp(\Ejj-\Ejpjp)}\gZmast q^{\pm(\Ejj-\Ejpjp)}\,,
\makebox[9em][r]{(6.3a)}\cr
\makebox[3em]{}
\xi(\Xpj)\cdot \gZmast & = &
-q^{-1}\Fm\Ejjp (\Em+\gZmast) \cr
& & \times [\Fm+q^{(\Ejj-\Ejpjp)/2}(\Em+\gZmast)]^{-1} \cr
& & - \Fm q^{-(\Ejj-\Ejpjp)/2}(\Em+\gZmast) \cr
& & \times [\Fm+q^{-(\Ejj-\Ejpjp)/2}(\Em+\gZmast)]^{-1} \cr
& & \times [\Fm-q^{-1}\Ejjp(\Em+\gZmast)] \cr
& & \times [\Fm+q^{(\Ejj-\Ejpjp)/2}(\Em+\gZmast)]^{-1} \cr
& = & -q^{-1}\Fm\Ejjp\gZmast+\gZmast\Ejjp\Em \cr
& & +\delta_{jm}q^{-1/2}q^{E_{m+1,m+1}/2}\gZmast E_{m,m+1}\gZmast
q^{-E_{mm}/2}\,,
\,\makebox[3em][r]{(6.3b)}\cr
\makebox[3em]{}
\xi(\Xmj)\cdot \gZmast & = &
-q\Fm\Ejpj (\Em+\gZmast) \cr
& & \times [\Fm+q^{(\Ejj-\Ejpjp)/2}(\Em+\gZmast)]^{-1} \cr
& & - \Fm q^{-(\Ejj-\Ejpjp)/2}(\Em+\gZmast) \cr
& & \times [\Fm+q^{-(\Ejj-\Ejpjp)/2}(\Em+\gZmast)]^{-1} \cr
& & \times [\Fm-q\Ejpj(\Em+\gZmast)] \cr
& & \times [\Fm+q^{(\Ejj-\Ejpjp)/2}(\Em+\gZmast)]^{-1} \cr
& = & -\delta_{jm}q^{1/2}E_{m+1,m}-q\Ejpj\gZmast+\gZmast\Ejpj\,.
\makebox[4em][r]{(6.3c)}\cr
\end{eqnarray*}
All manipulations needed here are quite straightforward. Note, for example,
that
$$
\Em[\Fm+q^{\pm(\Ejj-\Ejpjp)/2}(\Em+\gZmast)]^{-1}=
\Em q^{\mp(\Ejj-\Ejpjp)/2}\,.
$$

The relations (6.3), expressed directly in terms of entries $\zmast_{st}$,
were already presented in Ref. \cite{QG}. Let us recall them
($1\le s\le m<t\le N$):
\begin{eqnarray*}
\makebox[3em]{}
\xi(\qpmHj)\cdot\zmast_{st} & = &
q^{\pm(\delta_{js}-\delta_{j+1,s}-\delta_{jt}+\delta_{j+1,t})}\zmast_{st}\,;
\cr
\makebox[3em]{}
\xi(\Xpj)\cdot\zmast_{st} & = &
\delta_{j+1,s}\zmast_{jt}\,,\makebox[8em]{}\quad\ for\ j<m\,, \cr
& = & q^{-(1+\delta_{sm}-\delta_{t,m+1})/2}\zmast_{mt}\zmast_{s,m+1}\,,
\makebox[2em]{}\quad\ \ j=m\,,
\makebox[5em][r]{(6.4)}\cr
& = & -q^{-1}\delta_{jt}\zmast_{s,j+1}\,,\!\makebox[9em]{}\quad\ j>m\,; \cr
\makebox[3em]{}
\xi(\Xmj)\cdot\zmast_{st} & = &
\delta_{js}\zmast_{j+1,t}\,,\makebox[9em]{}\quad for\ j<m\,, \cr
& = & -q^{1/2}\delta_{sm}\delta_{t,m+1}\,,\,\makebox[7em]{}\,\quad\ j=m\,, \cr
& = & -q\delta_{j+1,t}\zmast_{sj}\,,\makebox[8em]{}\,\quad\ j>m\,. \cr
\end{eqnarray*}
It is quite straightforward to transcribe these relations in terms of
generators $\zeta^\ast_{st}$ for
$\zeta_{st}=z^{(s)}_{st},\  1\le s<t\le N$ (cf. (3.20) or (3.21)).
Note also that
$$
z^{(s)}_{s-1,t}=\zeta_{s-1,t}-\zeta_{s-1,s}\zeta_{st}\,.
$$
The result is given in (6.9) below.

It remains to determine the action on the unit. Recalling that
$M=\Lambda^\ast\Lambda$ and $\Lambda=(\alpha_{jk})$ is upper triangular,
one derives immediately from (4.26) that ($\lambda_{N+1}:=0$)
$$
\alpha_{st}\cdot1=0\quad for\ s<t,\ \alpha^2_{ss}\cdot1=\lambda_s,\
\alpha^\ast_{s,s+1}\alpha_{ss}\cdot1=
(\lambda_s-\lambda_{s+1})\zsast_{s,s+1}\,.
\eqno(6.5)$$
We introduce the substitution
$$
\lambda_{j+1}/\lambda_j=:q^{-2\sigma_j},\ j=1,\dots,N-1.
\eqno(6.6)$$
In virtue of (2.9), the expressions in Chevalley generators follow easily and
are given in (6.7) below.

Now we are ready to reformulate Proposition 4.4. Here we employ the
quantum coordinate functions on the flag manifold, i.e., the generators
$\zeta_{st}$ rather than $\zm_{st}$.
\proclaim Proposition 6.1.
A left \Uh -- module structure on \Fah depending on $(N-1)$ scalar
parameters $\sigma_1,\dots,\sigma_{N-1}$ is defined unambiguously by the
relations ($[x]:=(q^x-q^{-x})/(q-q^{-1})$)
$$
\qpmHj\cdot1=q^{\mp\sigma_j},\
\Xpj\cdot1=-q^{-(1+\sigma_j)/2}[\sigma_j]\zeta^\ast_{j,j+1},\
\Xmj\cdot1=0,\quad 1\le j\le N-1\,,
\eqno(6.7)$$
and
$$
Y\cdot(\zeta^\ast_{st}f)=\big(\xi(Y_{(1)})\cdot\zeta^\ast_{st}\big)\,
Y_{(2)}\cdot f,\quad with\ \DeltaYa\,,
\eqno(6.8)$$
for $1\le s<t\le N,\ \forall f\in\Faha$ and $\forall Y\in\Uha$. Here $\xi$
is the action dual to the right dressing transformation and it is prescribed
on the generators as follows ($1\le s<t\le N$):
\begin{eqnarray*}
\makebox[2em]{}
\xi(\qpmHj)\cdot\zeta^\ast_{st} & = &
q^{\pm(\delta_{js}-\delta_{j+1,s}-\delta_{jt}+\delta_{j+1,t})}
\zeta^\ast_{st}\,, \cr
\makebox[2em]{}
\xi(\Xpj)\cdot\zeta^\ast_{st} & = &
\delta_{j,s-1}(\zeta^\ast_{s-1,t}-\zeta^\ast_{st}\zeta^\ast_{s-1,s})
+\delta_{js}q^{-1+\delta_{s+1,t}/2}\zeta^\ast_{st}\zeta^\ast_{s,s+1}
-q^{-1}\delta_{jt}\zeta^\ast_{s,t+1}\,, \cr
\makebox[2em]{}
\xi(\Xmj)\cdot\zeta^\ast_{st} & = &
-q^{1/2}\delta_{js}\delta_{s,t-1}-q(1-\delta_{s,t-1})\delta_{j,t-1}
\zeta^\ast_{s,t-1} \makebox[10em][r]{(6.9)} \cr
& = & -q^{1-\delta_{s+1,t}/2}\delta_{j,t-1}\zeta^\ast_{s,t-1}\quad
(with\ \zeta^\ast_{ss}:=1)\,. \cr
\end{eqnarray*}

The unit is a lowest weight vector for $\Xmj\cdot1=0,\ \forall j$.
According to (6.7), the corresponding lowest weight equals
$$
\lambda:=-\sum_{j=1}^{N-1}\sigma_j\,\omega_j\,,
\eqno(6.10)$$
where $\omega_j$'s are the fundamental weights for {\gt su}$(N)$ defined by
$\omega_j(H_k)=\delta_{jk}$. Let \Ml designate the cyclic submodule of \Fah
with the cyclic vector $1$. It is known \cite{Lusztig, Rosso} that \Ml
is determined unambiguously, up to equivalence, by the lowest weight
$\lambda$ and the relation between finite-dimensional irreducible modules
and lowest weights is the same in the deformed as well as in the non-deformed case. This observation implies
\proclaim Proposition 6.2.
The unit in \Fah is a lowest weight vector corresponding to the lowest
weight $\lambda$ determined by $\lambda(H_j)=-\sigma_j,\ \forall j$.
The cyclic submodule $\Mla:=\Uha\cdot1$ in \Fah is a finite-dimensional
irreducible \Uh -- module provided all parameters $\sigma_j$ are
non-negative integers.


\begin{thebibliography}{article}

\bibitem{T-T} E. Taft and J. Towber, ``Quantum deformation of flag
schemes and Grassmann schemes I. A q-deformation  of the shape--algebra for
$GL(n)$,'' J. Algebra {\bf 142}, No. 1 (1991).

\bibitem{Soibelman} Y. Soibelman, ``On quantum flag manifolds,'' Funct. Anal.
Appl. {\bf 26}, 225 (1992).

\bibitem{Podlesz} P. Podlesz , ``Quantum spheres ,''
Lett. Math. Phys. {\bf 14}, 193 (1987).

\bibitem{Sheu} A. J. L. Sheu, ``Quantization of the Poisson $SU(2)$ and its
Poisson homogeneous space -- the 2-sphere,''
Commun. Math. Phys. {\bf 135}, 217 (1991).

\bibitem{JMP} P. \v S\v tov\'\i \v cek, ``Quantum line bundles on $S^2$
and the method of orbits for $SU_q(2)$,'' J. Math. Phys. {\bf 34},
1606 (1993).

\bibitem{QG} P. \v S\v tov\'\i\v cek, ``Quantum Grassmann manifolds,''
Commun. Math. Phys. {\bf 158}, 135 (1993).

\bibitem{J-S-II} B. Jur\v co and P. \v S\v tov\'\i\v cek,
``Coherent states for quantum compact groups,'' preprint CERN-TH. 7201/94,
and Commun. Math. Phys. (in press).

\bibitem{Parshall-W} B. Parshall and J. Wang, {\em Quantum linear groups}
(Rhode Island, AMS 1991).

\bibitem{Japanes} H. Awata, M. Noumi and S. Odake, ``Heisenberg realization
for $U_q(sl_n)$ on the flag manifold,''
Lett. Math. Phys. {\bf 30}, 35 (1994).

\bibitem{Biedenharn-L} L. C. Biedenharn and M. H Lohe, ``An extension of the
Borel-Weil construction to the quantum group $U_q(n)$,'' Commun. Math. Phys.
{\bf 146}, 483 (1992).

\bibitem{Fiore} G. Fiore, ``Realization of $U_q(so(N))$ within the
Differential Algebra on $R_q^N$,'' Commun. Math. Phys.
{\bf 169}, 475 (1995).

\bibitem{Brzezinski-M} T. Brzezi\' nski and S. Majid, ``Quantum group gauge
theory on quantum spaces,'' Commun. Math. Phys. {\bf 157}, 501 (1993).

\bibitem{Schuppetal} P. Schupp, P. Watts and B. Zumino, ``Cartan calculus for
Hopf algebras and quantum groups,'' preprint LBL--34215.

\bibitem{Woronowicz} S. L. Woronowicz, ``Differential calculus on quantum
matrix pseudogroups (quantum groups),''
Commun. Math. Phys. {\bf 122}, 125 (1989).

\bibitem{Jurco} B. Jur\v co, ``Differential calculus on quantized simple Lie
groups,'' Lett. Math. Phys. {\bf 22}, 177 (1991).

\bibitem{Chari-P} V. Chari and A. Pressley, {\em A guide to quantum groups}
(Cambridge University Press 1994).

\bibitem{Jurco-Sch} B. Jur\v co and M. Schlieker, ``On Fock-space
representations of quantized enveloping algebras related to noncommutative
differential geometry,'' J. Math. Phys. {\bf 36}, 3814 (1995).

\bibitem{Bulgariens} H. Sazdjian, Y. S. Stanev and I. T. Todorov,
``SU$_3$ coherent state operators and invariant correlation functions and
their quantum group counterparts,'' J. Math. Phys. {\bf 36}, 2030 (1995).

\bibitem{Kirillov-R} A. N. Kirillov and N. Yu. Reshetikhin, ``$q$--Weyl
group and a multiplicative formula for universal $R$--matrices,''
Commun. Math. Phys. {\bf 134}, 421 (1990).

\bibitem{Levendorskii-S} S. Z. Levendorskii and Ya. S. Soibelman,
``Some application of quantum Weyl groups. The multiplicative formula for
universal $R$--matrix for simple Lie algebra,''
J. Geom. Phys. {\bf 7}(4), 1 (1990).

\bibitem{Khoros-Tolstoy} S. M. Khoroshkin and V. N. Tolstoy, ``Universal
$R$--matrix for quantized (super)algebras,''
Commun. Math. Phys. {\bf 141}, 599 (1991).

\bibitem{SemenovTS} M. A. Semenov--Tian--Shansky, ``Dressing transformation
and Poisson Lie group actions,'' Publ. RIMS, Kyoto University {\bf 21},
1237 (1985).

\bibitem{Lu-W} J. H. Lu  and A. Weinstein, ``Poisson Lie groups, dressing
transformations and Bruhat decompositions,'' J. Diff. Geom. {\bf 31},
501 (1990).

\bibitem{JALG} P. \v S\v tov\'\i \v cek, ``Antiholomorphic representations
for orthogonal and symplectic quantum groups,'' J. Algebra (1996) (to appear).

\bibitem{F-R-T} N. Yu. Reshetikhin, L. A. Takhtajan and L. D. Faddeev,
``Quantization of Lie groups and Lie algebras,''
Algebra i analiz {\bf 1}, 178 (1989) (in Russian) and
Leningrad Math. J. {\bf 1}, 193 (1990).

\bibitem{Drinfeld} V. G. Drinfeld, ``Quantum groups,'' In Proc.
ICM Berkley 1986, (AMS 1987), p.798.

\bibitem{Jimbo1} M. Jimbo, ``A q--difference analogue of U(g) and the
Yang--Baxter equation,'' Lett. Math. Phys. {\bf 10}, 63 (1985).

\bibitem{Jimbo2} M. Jimbo, ``Quantum R-matrix for the generalized
Toda system,'' Commun. Math. Phys. {\bf 102}, 537 (1986).

\bibitem{J-S-I} B. Jur\v co and P. \v S\v tov\'\i \v cek,
``Quantum dressing orbits on compact groups,''
Commun. Math. Phys. {\bf 152}, 97 (1993).

\bibitem{Lusztig} G. Lusztig, ``Quantum deformations of certain simple
modules over enveloping algebras,'' Adv. Math. {\bf 70}, 237 (1988).

\bibitem{Rosso} M. Rosso, ``Finite dimensional representations for the
quantum analog of the enveloping algebra of a complex simple Lie algebra,''
Commun. Math. Phys. {\bf 117}, 581 (1988).

\end{thebibliography}
\end{document}